\documentclass[a4paper]{article}
\usepackage[latin1]{inputenc}
\usepackage[colorlinks=true]{hyperref}
\usepackage[T1]{fontenc}
\usepackage[francais,english]{babel}
\usepackage{eqsection}
\usepackage{graphics,amsmath,amssymb,amsfonts}
\usepackage{latexsym}
\usepackage{epsfig}
\usepackage{color}
\usepackage{color}
\normalbaselines
\newtheorem{theorem}{Theorem}[section]
\newtheorem{lemma}[theorem]{Lemma}
\newtheorem{corollary}[theorem]{Corollary}
\newtheorem{proposition}[theorem]{Proposition}

\newtheorem{remark}[theorem]{Remark}
\newtheorem{example}[theorem]{Example}

\newtheorem{definition}[theorem]{Definition}

\newenvironment{Proof}{\removelastskip\par\medskip
\noindent{\em Proof.} \rm}{\penalty-20\null\hfill$\square$\par\medbreak}

\def\Hob{\mathcal H_0}
\def\Hsta{\mathcal H_\beta}

\font\amsmath=msbm10 at 11pt

\newcommand{\R}{{\mbox{\amsmath R}}}

\title{\bf Control and stabilization of degenerate wave equations\footnote{This work has been partly supported by the GDRE CONEDP, a research group funded by CNRS and INdAM.}}
\author{Fatiha Alabau-Boussouira\thanks{
Laboratoire Jacques-Louis Lions, UMR-CNRS 7598, 
 Universit\'e de Lorraine, (France).
}
\and 
Piermarco
Cannarsa\thanks{Dipartimento di Matematica, Universit\`a di Roma Tor Vergata,
Via della Ricerca Scientifica, 00133 Roma (Italy); e-mail:
$<$cannarsa@axp.mat.uniroma2.it$>$. This work was completed while this author was visiting the Institut Henri Poincar\'e and Institut des Hautes \'Etudes Scientifiques on a senior  CARMIN position. }
\and G\"unter Leugering\thanks{Department Mathematik, LS Angewandte Mathematik 2,
Friedrich-Alexander-Universit\"at Erlangen-N\"urnberg, 
D-91058 Erlangen 
(Germany)}\thanks{The research of this author was supported by DFG EC315 "Engineering of Advanced Materials"}}

\begin{document}
%\date{}
\maketitle
\begin{abstract}
We study a wave equation in one space dimension with a general diffusion coefficient which degenerates on part of the boundary. Degeneracy is measured by a real parameter $\mu_a>0$.
We establish observability inequalities for weakly (when $\mu_a \in [0,1[$) as well as strongly (when $\mu_a \in [1,2[$) degenerate  equations. We also prove a negative result when the diffusion coefficient degenerates too violently (i.e. when $\mu_a>2$) and the blow-up of the observability time when $\mu_a$ converges to $2$ from below. Thus, using the HUM method we deduce the exact controllability of the corresponding degenerate control problem when $\mu_a \in [0,2[$. We conclude the paper by studying the boundary stabilization of the degenerate linearly damped wave equation and show that a suitable boundary feedback stabilizes the system exponentially. We extend this stability analysis to the degenerate nonlinearly boundary damped wave equation, for an arbitrarily growing nonlinear feedback close to the origin. This analysis proves that the degeneracy does not affect the optimal energy decay rates at large time. We apply the optimal-weight convexity method of \cite{alaamo2005, alajde2010} together with the results of the previous section, to perform this stability analysis.
\end{abstract}


\medskip


\medskip

{\bf Abbreviated title}. Control of degenerate wave equations

%%%%%%%%%%%%%%%%%%%%%%%%%%%%%%%%%%%%%%%%%%%%%%%%%%%%%%%%%%%%%%%%
%%%%%%%%%%%%%%%%%%%%%%%%%%%%%%%%%%%%%%%%%%%%%%%%%%%%%%%%%%%%%%%%
\section{Introduction}
%%%%%%%%%%%%%%%%%%%%%%%%%%%%%%%%%%%%%%%%%%%%%%%%%%%%%%%%%%%%%%%%
Control and inverse problems for degenerate PDE's arise in many applications such as cloaking (building of devices that lead to invisibility properties from observation) \cite{greenleaf}, climatology \cite{Diaz93}, population genetics \cite{Gurtinetal87}, and vision \cite{Manfredini04}. Such a variety of  applications has given birth to challenging mathematical problems for degenerate PDE's. A common feature of these problems is that they involve operators with variable diffusion coefficients that are not uniformly elliptic in the space domain, even though they are in general uniformly elliptic in compact subsets of the space domain, provided that these subsets are at a positive distance from the degeneracy.  This degeneracy may occur either on a part of the boundary or on a sub-manifold of the space domain. 

\medskip

The loss of uniform ellipticity  rises new questions related to the well-posedness of the evolution equations in suitable functional spaces as well as new estimates for the underlying elliptic equations. Similarly, in the degenerate case, new tools are necessary for the analysis of observability/nonobservability as well as stabilization.

\medskip

 Control issues for degenerate parabolic equations have received a lot of attention in the last ten years or so (see, for instance, \cite{MR2106292, Cannarsa-V-M-ADE, Cannarsa-V-M-SIAM, MR2538102}, \cite{Ala-Can-Fra}, \cite{Martinez-Vancost-JEE-2006}, and \cite{Can-Fra-Roc_1, Can-Fra-Roc_2}). New Carleman estimates with adapted weight functions, compared to the usual ones for nondegenerate parabolic equations, have been used to derive observability inequalities for the corresponding dual problems. 
 
 \medskip
 
Although degenerate wave equations have received less attention so far, we do believe that time has now come for a complete analysis and deeper understanding of these problems. Therefore, the purpose of this paper is to study controllability and observability issues for degenerate wave equations of the form
\begin{equation}
\label{intro_DW}
u_{tt} - \big(a(x)u_x\big)_x=0  \quad \mbox{in}\;\; ]0,\infty[\times]0,1[,
\end{equation}
where $a$ is positive on $]0,1]$ but vanishes at zero. Moreover,  if stabilization is usually irrelevant in the parabolic case due to the intrinsic dissipation of diffusion models, it remains an important question for degenerate wave equations and will be addressed in this paper.

\medskip

The degeneracy of \eqref{intro_DW} at $x=0$ is measured by the parameter $\mu_a$ defined by 
\begin{equation}\label{intro:hp_a}
 \mu_a:=\sup_{0<x\leqslant 1}\dfrac{x|a'(x)|}{a(x)},
 \end{equation}
 and one says that \eqref{intro_DW} degenerates weakly if $\mu_a \in [0,1[$, strongly if $\mu_a>1$. 
Here we assume $\mu_a<2$ because,  like in the parabolic case, observability no longer holds true if $\mu_a \geqslant 2$ as we show in Section~\ref{se:failure} of this paper.
 
By determining suitable multipliers linked to the coefficient $\mu_a$ of the degeneracy and proving refined trace theorems, we prove  boundary observability inequalities for \eqref{intro_DW} in a sufficiently large time. This approach and tools are new in the context of degenerates wave equations, as far as we know.
 It is worth noting that, in problems involving cloaking which, obviously, is incompatible with observability, the degeneracy of the coefficients is quadratic (see \cite{greenleaf}). So, our results are consistent with such a property.  Moreover, we also study the behavior of the controllability (or observability) time as $\mu_a$ converges to $2$,  appealing to Bessel functions to show that such a time blows up as
$\mu_a$ approaches to $2$ from below. 

For a certain class of weakly degenerate wave equations, an interesting result with $x=0$ as observation region was obtained by Gueye~\cite{gueye} by using the explicit description of the spectrum of the corresponding elliptic operator to treat the related moment problem. As a consequence, an exact controllability result with Dirichlet boundary controls located at the degeneracy point was deduced for degenerate wave equations and then extended to degenerate parabolic equations, giving a first answer to a question that had been open for quite some time. The viewpoint of this paper differs from the one of \cite{gueye}. Indeed, we allow coefficients to degenerate either weakly or strongly on the boundary and we  obtain observability or controllability from the nondegenerate part of the boundary. Moreover, we employ direct techniques such as multipliers and sharp trace results. 

Finally, we devote a substantial part of the paper to the study of boundary stabilization for \eqref{intro_DW} when $\mu_a \in [0,1[$. We consider both the linear feedback
\begin{equation}\label{intro_lf}
u_t(t,1)+u_x(t,1)+\beta u(t,1)=0,
\end{equation}   
and the nonlinear damping 
\begin{equation}\label{intro_nlf}
\rho(u_t(t,1))+u_x(t,1)+\beta u(t,1)=0,
\end{equation}
where $\beta > 0$, and $\rho$ is a nondecreasing function such that $\rho(0)=0$. Thanks to the dominant energy approach together with suitable elliptic estimates, we prove that \eqref{intro_lf}  stabilizes exponentially the corresponding solution of the degenerate wave equation. For the nonlinear feedback \eqref{intro_nlf}, we use the optimal-weight convexity method of \cite{alaamo2005, alajde2010} to establish a quasi-optimal energy decay rate using the multipliers we have determined in the linear case. We also discuss several explicit examples of decay corresponding to different feedbacks. 
We recall that, for finite dimensional models, the optimality of the decay rates provided by the optimal-weight convexity method is proved  in \cite{alajde2010}. Moreover, our results show that, under the action of a nonlinear boundary damping, degenerate wave equations enjoy the same stability properties as the corresponding nondegenerate equations, in the sense that both models have the same decay rates of the energy. 

We would like to point out that one can reformulate all the above results on nonlinear stabilization by replacing integral inequalities with a Lyapunov function technique. As we explain in Remark~\ref{re:Lyapunov} below, such an operation is essentially of no use.
\medskip

The paper is organized as follows. In section 2, we introduce our notations, define the degeneracy parameter $\mu_a$, functional spaces and assumptions. We also prove Poincar\'e's type inequalities and some key trace results for functions in  weighted Sobolev spaces. In section 3, we consider the dual problem, prove well-posedness, and prove the direct inequality as well as
the boundary observability property for $ \mu_a \in [0,2[$. We prove non-observability for $\mu_a>2$ and the blow-up of the observability time when $\mu_a$ converges to $2$ from below.
We conclude this section by proving exact boundary controllability for the controlled system when  $ \mu_a \in [0,2[$. We consider the boundary stabilization problem in section 4 and
prove its well-posedness, together with its exponential stability. We extend this stability analysis to the nonlinear boundary stabilization problem in section 5.

%%%%%%%%%%%%%%%%%%%%%%%%%%%%%%%%%%%%%%%%%%%%%%%%%%%%%%%%%%%%%%%%
%%%%%%%%%%%%%%%%%%%%%%%%%%%%%%%%%%%%%%%%%%%%%%%%%%%%%%%%%%%%%%%%
\section{Assumptions and preliminaries}
%%%%%%%%%%%%%%%%%%%%%%%%%%%%%%%%%%%%%%%%%%%%%%%%%%%%%%%%%%%%%%%%
\subsection{Assumptions}
Let $a\in\mathcal C([0,1])\cap \mathcal C^1(]0,1])$ be a function satisfying the following assumptions:
\begin{equation}\label{eq:hp_a}
\begin{cases}
(i)&  a(x)>0\;\;\forall x\in ]0,1]\,,\;\;
a(0)=0\,,
\\
(ii)& \mu_a:=\sup_{0<x\leqslant 1}\dfrac{x|a'(x)|}{a(x)}<2 \,,\;\;\mbox{and}
\\
(iii)&
a\in \mathcal C^{[\mu_a]}([0,1]),
\end{cases}
\end{equation}
where $[\cdot]$ stands for the integer part.
\begin{remark}\label{re:a}\rm Assumption~\eqref{eq:hp_a} subsumes similar hypotheses that were formulated to treat degenerate parabolic equations (see, for instance, \cite{Ala-Can-Fra, Cannarsa-V-M-SIAM, Martinez-Vancost-JEE-2006}). We list below some simple consequences of \eqref{eq:hp_a}.
\begin{enumerate}
\item By integrating the inequality  
\begin{equation*}
%\label{eq:}
sa'(s)\leqslant \mu_a a(s)\qquad\forall\,s\in]0,1]
\end{equation*}
over $[x,1]$ we obtain 
\begin{equation}\label{eq:a}
a(x)\geqslant a(1) x^{\mu_a}\quad\forall x\in [0,1]\,.
\end{equation}
Consequently,  $1/a\in L^1(0,1)$ when $\mu_a\in [0,1[$.  
\item Observe that condition \eqref{eq:hp_a} $(iii)$ is equivalent to require that $a\in \mathcal C^{1}([0,1])$ when $\mu_a\in [1,2[$ (no extra assumption is imposed when $\mu_a\in [0,1[$).
In this case of strong degeneracy, we have that $1/a\notin L^1(0,1)$. Indeed, since $a\in \mathcal C^{1}([0,1])$, we have that  
\begin{equation*}
%\label{eq:}
\dfrac{a(x)}x<1+|a'(0)|
\end{equation*}
in some neighborhood of $0$. So, $1/a\notin L^1(0,1)$. 
\end{enumerate}
\end{remark}
\subsection{Function spaces}
We now introduce some weighted Sobolev spaces that are naturally associated with degenerate operators, see \cite{Cannarsa-V-M-SIAM}. We denote by $V^1_a(0,1)$ the space of all functions $u\in L^2(0,1)$ such that
\begin{equation}\label{eq:H1a}
\begin{cases}
(i)&
u\;\;\mbox{is locally absolutely continuous in}\;]0,1],
%\\
%(ii)&
%u(1)=0, 
\;\mbox{and}
\\
(iii)&
\sqrt{a}u_x \in L^2(0,1).
\end{cases}
\end{equation}
It is easy to see that $V^1_a(0,1)$ is an Hilbert space with the scalar product
\begin{equation*}
%\label{eq:}
\langle u,v\rangle_{1,a}= \int_0^1\big(a(x)u'(x)v'(x)+u(x)v(x)\big)dx\,,\qquad\forall\,u,v\in V^1_a(0,1)
\end{equation*}
and associated norm
\begin{equation*}
%\label{eq:}
\|u\|_{1,a}= \Big\{\int_0^1\big(a(x)|u'(x)|^2+|u(x)|^2\big)dx\Big\}^{\frac12}\,,\qquad\forall\,u\in V^1_a(0,1)\,.
\end{equation*}
Let us also set
\begin{equation*}
%\label{eq:}
|u|_{1,a}=\Big\{\int_0^1a(x)\,|u'(x)|^2dx\Big\}^{\frac12}\qquad\forall\,u\in V^1_a(0,1).
\end{equation*}
Actually, $|\cdot|_{1,a}$ is an equivalent norm on the closed subspace of $V^1_{a,0}(0,1)$ defined as
\begin{equation*}
%\label{eq:}
V^1_{a,0} (0,1)=\big\{u \in V^1_a(0,1)~:~u(1)=0\big\}.
\end{equation*}

This fact is a simple consequence of the following version of Poincar\'e's inequality.
\begin{proposition}\label{pr:poincare}
 Assume \eqref{eq:hp_a}.
 Then 
\begin{equation}\label{eq:poincare}
\|u\|^2_{L^2(0,1)}\leqslant C_a\; |u|^2_{1,a}\qquad\forall\,u\in V^1_{a,0} (0,1),
\end{equation}
where
\begin{equation}\label{eq:cost_pcr}
C_a=\frac{1}{a(1)} \min \Big\{ 4, \dfrac{1}{ 2-\mu_a}\Big\}\,.
\end{equation}
\end{proposition}
\begin{Proof} Let $u\in V^1_{a,0} (0,1)$. We will prove two different bounds for $\|u\|^2_{L^2(0,1)}$ in terms of $|u|^2_{1,a}$. The conclusion \eqref{eq:poincare} will follow by taking the minimum of the two corresponding constants.

First, we use a direct argument. For any $x\in ]0,1]$ we have that
\begin{equation*}
%\label{eq:}
|u(x)|=\Big|\int_x^1u'(s)ds\Big|\leqslant |u|_{1,a}\, \Big\{\int_x^1\dfrac{ds}{a(s)}\Big\}^{\frac12}.
\end{equation*}
Therefore, by Fubini's theorem,
\begin{eqnarray}\label{eq:fubini}
\int_0^1|u(x)|^2dx\leqslant |u|^2_{1,a}\, \int_0^1dx\int_x^1\dfrac{ds}{a(s)} =|u|^2_{1,a}\, \int_0^1\dfrac{ds}{a(s)} \int_0^sdx=
|u|^2_{1,a}\, \int_0^1\dfrac{s}{a(s)} ds\,.
\end{eqnarray}
By \eqref{eq:a} we deduce that 
\begin{eqnarray*}
\int_0^1\dfrac{s}{a(s)} ds\leqslant \frac 1{a(1)}\int_0^1s^{1-\mu_a } ds\,.
\end{eqnarray*}
Together with \eqref{eq:fubini}, the above inequality yields the first bound we mentioned above, that is,
 \begin{equation}\label{eq:poincare1}
\|u\|^2_{L^2(0,1)}\leqslant \dfrac{\; |u|^2_{1,a}}{a(1)(2-\mu_a)}\qquad\forall\,u\in V^1_{a,0} (0,1).
\end{equation}

Next, as an alternative proof, we adapt a reasoning that can be used to prove Hardy's inequality. Observe that, for all $x\in]0,1[$,
\begin{eqnarray*}
0&\leqslant&\int_x^1\Big(su'(s)+\frac12 u(s)\Big)^2ds
\\
&=&\int_x^1\Big(s^2|u'(s)|^2+\frac14 |u(s)|^2+s\,u(s)\,u'(s)\Big) ds
\\
&=&\int_x^1\Big(s^2|u'(s)|^2-\frac14 |u(s)|^2 \Big) ds-\frac12 x\,|u(x)|^2.
\end{eqnarray*}
Therefore, taking the limit as $x\downarrow 0$, by \eqref{eq:a} we obtain the announced second  bound:
\begin{equation}
\label{eq:poincare2}
\int_0^1|u(s)|^2 ds\leqslant 4\int_0^1 s^2|u'(s)|^2ds\leqslant 4\int_0^1 s^{\mu_a}|u'(s)|^2ds\leqslant \frac4{a(1)}\int_0^1 a(s)|u'(s)|^2ds\qquad\forall\,u\in V^1_{a,0} (0,1).
\end{equation}
The conclusion follows from \eqref{eq:poincare1} and \eqref{eq:poincare2}.
\end{Proof}

\begin{example}\label{exa:power}\rm The following are examples of functions $a$ satisfying assumption \eqref{eq:hp_a}. 
\begin{enumerate}
\item Let  $\theta\in]0,2[$ be given. Define
\begin{equation}\label{eq:power}
a(x)=x^\theta\qquad \forall\,x\in[0,1].
\end{equation}
In this case, we have 
\begin{equation}
\label{eq:poincare_power}
\|u\|^2_{L^2(0,1)}\leqslant \min \Big\{ 4, \dfrac{1}{ 2-\theta}\Big\}\; |u|^2_{1,a}\qquad\forall\,u\in V^1_{a,0} (0,1).
\end{equation}
\item Let  $\theta\in]0,2[$ be given and let $\alpha\in]0,1-\theta/2[$. Then the function
\begin{equation}\label{eq:power+}
a(x)=
\begin{cases}
\hspace{.cm}
x^\theta\big(1+\sin^2(\log x^\alpha)\big)
&\forall\,x\in]0,1]
\\
\hspace{.cm}
0& x=0
\end{cases} 
\end{equation}
satisfies \eqref{eq:hp_a}. Indeed, 
\begin{equation*}
%\label{eq:}
a'(x)=\theta x^{\theta-1}\big(1+\sin^2(\log x^\alpha)\big)+2\alpha x^{\theta-1} \sin(\log x^\alpha) \cos(\log x^\alpha)\qquad \forall\,x\in]0,1]\,,
\end{equation*}
so that $\mu_a\leqslant  \theta+2\alpha<2$.
Notice that \eqref{eq:poincare_power} is still vaild for this weight function because
\begin{equation*}
%\label{eq:}
a(x)\geqslant  x^{\theta}\quad\forall x\in [0,1]\,,
\end{equation*}
which is what is really needed for the prooof of Proposition~\ref{pr:poincare}.
\end{enumerate}
\end{example}

\begin{remark}\label{re:poincare}\rm
We do not expect the constant $C_a$ in \eqref{eq:cost_pcr} to be optimal. For instance,  for the weight $a$ in \eqref{eq:power}, we have that
\begin{equation*}
%\label{eq:}
\min\Big\{ 4, \dfrac{1}{ 2-\theta}\Big\}\to \frac 12\quad\mbox{as}\quad \theta\downarrow 0\,,
\end{equation*}
which is strictly greater than the minimal constant in the case $\theta=0$, that is, $(2/\pi)^2$. On the other hand,  \eqref{eq:cost_pcr} shows that  $C_a$ does not blow up as $\mu_a\uparrow 2$ because it is  bounded above by $4$.
\end{remark}

Next, we define
$$
V^2_a(0,1)=\big\{u \in V^1_a(0,1)~:~au' \in H^1(0,1)\big\},
$$
where $H^1(0,1)$ denotes the classical Sobolev space of all functions $u\in L^2(0,1)$ such that $u'\in L^2(0,1)$. Notice that, if $u\in V^2_a(0,1)$, then $au'$ is continuous on $[0,1]$.

We collect below useful properties of the above functional spaces. Some of the following results are known, others are new. We prove all of them for completeness.
\begin{proposition}\label{pr:sobolev}
Assume \eqref{eq:hp_a}. Then the following properties hold true.
\begin{description}
\item[(I)] For every $u\in V^1_a(0,1)$ 
\begin{equation}\label{eq:H1a_0}
\lim_{x\downarrow 0}x\,u^2(x)=0,
\end{equation}
\begin{equation}\label{eq:trace}
u^2(1)\leqslant \max\Big\{2,\frac1{a(1)}\Big\}\|u\|^2_{1,a}\,.
\end{equation}

Moreover, if $\mu_a\in[0,1[$, then $u$ is  absolutely continuous in $[0,1]$.
\item[(II)] For every $u\in V^2_a(0,1)$
\begin{eqnarray}
\label{eq:H2a_1}
&\lim_{x\downarrow 0} &x\,a(x)\,u'(x)^2=0\,.
\end{eqnarray}
For all $u\in V^2_a(0,1)$ and $\phi\in V^1_a(0,1)$
\begin{equation}
\label{eq:H2a_2}
\lim_{x\downarrow 0} a(x)\,\phi(x)\,u'(x)=0\,,
\end{equation}
assuming, in addition, $\phi(0)=0$ when $\mu_a\in [0,1[$.
\item[(III)]
If $\mu_a\in[1,2[$, then for every $u\in V^2_a(0,1)$
\begin{equation}\label{eq:H2a_3}
\lim_{x\downarrow 0} 
a(x)\,u'(x)=0\,.
\end{equation}
\end{description}
\end{proposition}
\begin{Proof}
{\bf (I)} Let $u\in V^1_a(0,1)$. We will show that \begin{equation*}
%\label{eq:}
v(x):=
\begin{cases}
\hspace{.cm}x\,u^2(x)
&0<x\leqslant 1
\\
\hspace{.cm}
0 & x=0
\end{cases}
\end{equation*}
is  continuous on $[0,1]$. Indeed, $v$ is locally absolutely continuous in $]0,1]$ and
\begin{eqnarray*}
v'(x)=u^2(x)+2\,x\,u'(x)\,u(x) \qquad \mbox{a.e. in }[0,1].
\end{eqnarray*}
Now, the above right-hand side is in $L^1(0,1)$ because $u\in L^2(0,1)$ and, thanks to \eqref{eq:a},
\begin{equation}
\label{eq:late}
\int_0^1x^2u'(x)^2dx\leqslant \int_0^1x^{\mu_a} |u'(x)|^2dx\leqslant \dfrac1{a(1)}\int_0^1a(x) |u'(x)|^2dx\,.
\end{equation}
Then, the limit $\lim_{x\downarrow 0} v(x)=:L$ does exist and must vanish for otherwise $u^2(x)\sim L/x$ (near zero) would not be summable.
\eqref{eq:H1a_0} is thus proved.

Next, we have that
\begin{equation*}
%\label{eq:}
u^2(1)=v(1)=\int_0^1\big(u^2(x)+2\,x\,u'(x)\,u(x)\big)dx \leqslant 2\int_0^1u^2(x)dx+\int_0^1x^2|u'(x)|^2dx
\end{equation*}
which, in turn, yields 
\begin{equation*}
%\label{eq:}
u^2(1) \leqslant 2\int_0^1u^2(x)dx+\dfrac1{a(1)}\int_0^1a(x) |u'(x)|^2dx
\end{equation*}
in view of \eqref{eq:late}.

Now, suppose, in addition, that $\mu_a\in[0,1[$. Then
\begin{equation*}
%\label{eq:}
u'(x) = \dfrac{1}{\sqrt{a(x)}}\,\sqrt{a(x)}u'(x)\qquad\forall x\in ]0,1]
\end{equation*}
is summable over $(0,1)$ thanks to Remark~\ref{re:a}-1 and \eqref{eq:H1a}$(iii)$. So, $u$ is  absolutely continuous in $[0,1]$. 

{\bf (II)} Let $u\in V^2_a(0,1)$. We claim that 
\begin{equation*}
%\label{eq:}
v(x):=
\begin{cases}
\hspace{.cm}
x\,a(x)\,u'(x)^2
&0<x\leqslant 1
\\
\hspace{.cm}
0& x=0
\end{cases}
\end{equation*}
is  continuous on $[0,1]$. Indeed, $v$ is locally absolutely continuous in $]0,1]$ and
\begin{eqnarray*}
v'(x)&=&a(x)u'(x)^2+x\,a'(x)\,u'(x)^2+2\,x\,a(x)\,u'(x)\,u''(x)
\\
&=&a(x)u'(x)^2+2\,x\,u'(x)\big(a(x)u'(x)\big)'-x\,a'(x)\,u'(x)^2 \qquad \mbox{a.e. in }[0,1].
\end{eqnarray*}
Now, observe that the first term in the right-hand side above is summable over $[0,1]$ in view of \eqref{eq:H1a} $(iii)$, and the same is true for second one   because, by \eqref{eq:a},
\begin{equation*}
%\label{eq:}
x|u'(x)|\leqslant x^{\mu_a/2}|u'(x)|\leqslant \sqrt{\dfrac{a(x)}{a(1)}}\,|u'(x)|\qquad \forall\,x\in ]0,1]\,.
\end{equation*}
As for the third term, owing to \eqref{eq:hp_a} $(ii)$, 
\begin{eqnarray*}
 x\,|a'(x)|\,u'(x)^2\leqslant \mu_a \;a(x) \,u'(x)^2\qquad \forall\,x\in ]0,1]
\end{eqnarray*}
and the above right-hand side is summable in view of  \eqref{eq:a} $(iii)$. Then, $\lim_{x\downarrow 0} v(x)=:L$ does exist and must vanish, for otherwise $a(x)u^2(x)\sim L/x$ (near zero) would not be summable. This concludes the proof of 
\eqref{eq:H2a_1}.

{\bf (III)} Next, we  proceed to prove \eqref{eq:H2a_3} noting that 
 $\lim_{x\downarrow 0} a(x)\,u'(x)=:L$ exists because  $u\in V^2_a(0,1)$ and  must vanish, for otherwise $a(x)u'(x)^2\sim L^2/a(x)$ (near zero) would not be summable in view of Remark~\ref{re:a}-3.

 Finally, in order to show \eqref{eq:H2a_2}, we begin by proving that the function
 \begin{equation*}
%\label{eq:}
w(x):=
\begin{cases}
\hspace{.cm}
a(x)\,\phi(x)\,u'(x)
&0<x\leqslant 1
\\
\hspace{.cm}
0& x=0
\end{cases}
\end{equation*}
is  continuous on $[0,1]$. This follows by the arguments as above, because
\begin{eqnarray*}
w'(x)=a(x)\phi'(x)u'(x)+ \phi(x)\big(a(x)u'(x)\big)'
%\qquad\mbox{a.e. in}\;\;[0,1]\,
\end{eqnarray*}
 is summable over $[0,1]$. Therefore, the limit $\lim_{x\downarrow 0} w(x)=:L$  exists and $a(x)|\phi(x)\,u'(x)|\sim |L| $ near $0$. We now have to distiguish two cases. If $\mu_a\in [0,1[$ and $\phi(0)=0$, then the conclusion is immediate. If, on the other hand, $\mu_a\in [1,2[$, then, owing to \eqref{eq:H2a_3},
\begin{eqnarray*}
 a(x)|u'(x)|=\Big|\int_0^x\big(a(x)u'(x)\big)'dx\Big|\leqslant \sqrt{x}\,\big\|(au')'\big\|_{L^2(0,1)}\qquad\forall x\in]0,1]\,.
\end{eqnarray*}
Now, if $L\neq 0$, then, in a neighborhood of $0$, 
\begin{eqnarray*}
\frac{|L|}2\leqslant a(x)|\phi(x)\,u'(x)| \leqslant  \sqrt{x}\,\big\|(au')'\big\|_{L^2(0,1)}|\phi(x)|\,,\qquad\forall x\in]0,1]\,
\end{eqnarray*}
in contrast to the fact that $\phi\in L^2(0,1)$.
\end{Proof}

%%%%%%%%%%%%%%%%%%%%%%%%%%%%%%%%%%%%%%%%%%%%%%%%%%%%%%%%%%%%%%%%
%%%%%%%%%%%%%%%%%%%%%%%%%%%%%%%%%%%%%%%%%%%%%%%%%%%%%%%%%%%%%%%%
\section{Observability}
%%%%%%%%%%%%%%%%%%%%%%%%%%%%%%%%%%%%%%%%%%%%%%%%%%%%%%%%%%%%%%%%
Given $a$ satisfying assumptions \eqref{eq:hp_a}, let $\mu_a\in [0,2[$ be the constant in assumption $(ii)$. Consider the degenerate wave equation
\begin{equation}
\label{DW}
u_{tt} - \big(a(x)u_x\big)_x=0  \quad \mbox{in}\;\; ]0,\infty[\times]0,1[
\end{equation}
with
\begin{equation}
\label{b+i}
\begin{cases}
\mbox{boundary conditions}\;\;
u(t,1)=0 \;\;\mbox{and}\;\;\begin{cases}
\hspace{.cm}
u(t,0)=0
&
\mbox{if}\;\;\mu_a\in[0,1[
\\
\hspace{.cm}
\lim_{x\downarrow 0} 
a(x)\,u_x(t,x)=0
&
\mbox{if}\;\;\mu_a\in[1,2[
\end{cases}
&0<t<\infty
\\
\mbox{initial conditions}\;\;
\begin{cases}
\hspace{.cm}u(0,x)=u_0(x)
\\
u_t(0,x)= u_1(x)
\hspace{.cm}
\end{cases}
& x\in]0,1[.
\end{cases}
\end{equation}
We recall that, since  equation \eqref{DW} is  degenerate, different boundary conditions have to be imposed at $x=0$ depending on whether we are interested in:
\begin{itemize}
\item the {\em weakly degenerate} case $\mu_a\in[0,1[$, where, in view of Proposition~\ref{pr:sobolev}-(I), we have that the Dirichlet boundary condition $u(t,0)=0$ makes sense for any solution, and
\item the {\em strongly degenerate} case $\mu_a\in[1,2[$, where , in view of Proposition~\ref{pr:sobolev}-(II), we have that the Neumann boundary condition $\lim_{x\downarrow 0} 
a(x)\,u_x(t,x)=0$ is automatically satisfied by any classical solution. 
\end{itemize}
In order to express the  above boundary conditions in functional settings, we  define $H^1_{a}(0,1)$ to be the closed subspace of $V^1_{a,0} (0,1)$ which consists of all 
$u\in V^1_{a,0} (0,1)$ satisfying $u(0)=0$ when $\mu_a\in [0,1[$. We  also set 
$$
H^2_a(0,1)= V^2_a(0,1) \cap H^1_a(0,1).
$$
Observe that all functions $u\in H^2_a(0,1)$ satisfy  homogeneous boundary conditions at both $x=0$ and $x=1$. Such conditions are of Dirichlet type in the weakly degenerate case, whereas they are of Neumann/Dirichlet type at $x=0$ and $x=1$, respectively, when $\mu_a\in [1,2[$.

\subsection{Well-posedness}
Let us recall the typical abstract set-up of semigroup theory which provides weak and classical notions of solutions for problem \eqref{DW}-\eqref{b+i}. Consider the Hilbert space $\Hob=H^1_a(0,1)\times L^2(0,1)$ with the scalar product
\begin{equation*}
%\label{eq:}
\big\langle (u,v),(\widetilde u, \widetilde v)\big\rangle=
\int_0^1\big(v(x)\widetilde v(x)+a(x)u'(x)\widetilde u'(x)\big)dx
\qquad\forall\,(u,v),(\widetilde u,\widetilde v)\in \Hob\,.
\end{equation*}
Arguing as for the classical wave equation (see, for instance, \cite{tanabe}) one can show that the unbounded operator $A:D(A)\subset \Hob\to \Hob$ defined by
\begin{equation*}
%\label{eq:}
\begin{cases}
\hspace{.cm}
D(A)= H^2_a(0,1)\times H^1_a(0,1)
&
\\
A(u,v)=\big(v, (au')'\big)\qquad\forall (u,v)\in D(A)
\hspace{.cm}
\end{cases}
\end{equation*}
 is maximal dissipative  on $\Hob$. Therefore, $A$ is the generator of a contraction semigroup in $\Hob$, denoted by $e^{tA}$. For any $U_0:=(u_0,v_0)\in \Hob$, $U(t):= e^{tA}U_0$  gives the  so-called mild solution of the Cauchy problem
\begin{equation*}
%\label{eq:Cauchy}
\begin{cases}
\hspace{.cm}
U'(t)=AU(t)
&(t\geqslant 0)
\\
\hspace{.cm}
U(0)=U_0.
\end{cases}
\end{equation*}
When $U_0\in D(A)$, the above solution is classical  in the sense that $U\in \mathcal C^1\big([0,\infty[;\Hob\big)\cap \mathcal C\big([0,\infty[;D(A)\big)$ and the equation holds on $[0,\infty[$. 

In view of the above considerations, given $(u_0,u_1)\in H^1_a(0,1)\times L^2(0,1)$, we say that the function
\begin{equation*}
%\label{eq:}
u\in \mathcal C^1\big([0,\infty[;L^2(0,1)\big)\cap \mathcal C\big([0,\infty[;H^1_a(0,1)\big)
\end{equation*}
is the {\em mild} (or {\em weak}) {\em solution} of problem \eqref{DW}-\eqref{b+i} if  $\big(u(t),v(t)\big)=e^{tA}(u_0,v_0)$ for all $t\geqslant 0$. 
By the aforementioned regularity result for $e^{tA}$, if $(u_0,u_1)\in H^2_a(0,1)\times H^1_a(0,1)$, then $u$ is the {\em classical solution} of
\eqref{DW}-\eqref{b+i} meaning that
\begin{equation*}
%\label{eq:}
u\in \mathcal C^2\big([0,\infty[;L^2(0,1)\big)\cap \mathcal C^1\big([0,\infty[;H^1_a(0,1)\big)\cap \mathcal C\big([0,\infty[;H^2_a(0,1)\big)
\end{equation*}
and \eqref{DW}  is satisfied for all $t\in [0,\infty[$ and a.e. $x\in [0,1]$.

The energy of a mild solution $u$ of \eqref{DW} is the continuous function defined by
\begin{equation}\label{defenergy}
E_u(t)=\frac{1}{2} \int_0^1\big\{ u_t^2(t,x) + a(x) u_x^2(t,x)\big\}dx \qquad\forall t\geqslant 0\,.
\end{equation}
\begin{proposition}\label{pr:energy}
Assume \eqref{eq:hp_a} and let $u$ be the mild solution of \eqref{DW}-\eqref{b+i}. Then 
\begin{equation}\label{conserv}
E(t)=E(0) \quad \forall \ t \geqslant 0 \,.
\end{equation}
\end{proposition}
\begin{Proof} Suppose, first, that $u$ is a classical solution of \eqref{DW}. Then, multiplying the equation by $u_t$ and integrating by parts we obtain
\begin{eqnarray*}
0&=&\int_0^1u_t(t,x)\big\{u_{tt}(t,x)-\big(a(x)u_x(t,x)\big)_x\big\}dx
\\
&=&\underbrace{\int_0^1\big\{u_t(t,x)u_{tt}(t,x)+ a(x)u_x(t,x)u_{tx}(t,x) \big\} dx}_{=\frac d{dt}E_u(t)}
-\big[a(x)u_t(t,x)u_x(t,x)\big]_{x=0}^{x=1}.
\end{eqnarray*}
By noting that the boundary terms vanish because of the boundary conditions in both the weakly and strongly degenerate cases, we conclude that the energy of $u$ is constant. The same conclusion can be extended to any mild solution by an approximation argument.
\end{Proof}
 \subsection{Boundary observability}
\begin{lemma}\label{le:1}
For any mild solution  $u$ of \eqref{DW}  we have that $u_x(\cdot,1)\in L^2(0,T)$ for every $T\geqslant 0$ and
\begin{equation}\label{eq:le0}
a(1)\int_0^T u_x^2(t,1)\,dt\leqslant \Big(6T+\dfrac{1}{\min\{1, a(1)\}}\Big)E_u(0).
\end{equation}
Moreover,
\begin{equation}\label{eq:le1}
a(1)\int_0^T u_x^2(t,1)\,dt= \int_0^T\!\!\!\int_0^1\Big\{u_t^2(t,x) + \big(a(x)-xa'(x)\big)u_x^2(t,x)\Big\} dtdx+ 2\Big[\int_0^1 xu_x(t,x) u_t(t,x) dx\Big]_{t=0}^{t=T}.
\end{equation}
\end{lemma}
\begin{Proof}
Suppose first $(u_0,u_1)\in H^2_a(0,1)\times H^1_a(0,1)$ so that $u$ is a classical solution of \eqref{DW}. Then, by multiplying  equation \eqref{DW} by $xu_x$ and integrating over $]0,T[\times]0,1[$ we obtain
\begin{eqnarray}
\nonumber
0&=&\int_0^T\!\!\!\int_0^1 xu_x(t,x)\Big(u_{tt}(t,x) - \big(a(x)u_x(t,x)\big)_x\Big)dx\,dt
\\
\nonumber
&=&\Big[\int_0^1 xu_x(t,x) u_t(t,x) dx\Big]_{t=0}^{t=T}-\int_0^T\!\!\!\int_0^1xu_{tx}(t,x)\,u_{t}(t,x)\,dx\,dt
\\
\nonumber
& &\qquad-\int_0^T\!\!\!\int_0^1 \Big( x\,a'(x)\,u_x^2(t,x)+x\,a(x)\,u_x(t,x)u_{xx}(t,x)\Big)dx\,dt
\\
\nonumber
&=&\Big[\int_0^1 xu_x(t,x) u_t(t,x) dx\Big]_{t=0}^{t=T}-\int_0^T\!\!\!\int_0^1x\,a'(x)\,u_x^2(t,x)\,dx\,dt
\\
\label{eq:le1_id1}
& &\qquad-\int_0^T\!\!\!\int_0^1 \Big\{ x\Big(\dfrac{u^2_t(t,x)}2\Big)_x+x\,a(x)\Big(\dfrac{u^2_x(t,x)}2\Big)_x\Big\}dx\,dt
\end{eqnarray}
We  proceed to integrate by parts the last two terms above. We obtain
\begin{equation}
\label{eq:l1_m1}
\int_0^T\!\!\!\int_0^1 x\Big(\dfrac{u^2_t(t,x)}2\Big)_x\,dx\,dt=-\dfrac{1}2\int_0^T\!\!\!\int_0^1 u^2_t(t,x)\,dx\,dt
\end{equation}
because $xu^2_t(t,x)$ vanishes at $x=1$ and, owing to \eqref{eq:H1a_0}, also at $x=0$. Moreover, on account of \eqref{eq:H2a_1} we have
\begin{equation}
\label{eq:l1_e2}
\int_0^T\!\!\!\int_0^1x\,a(x)\Big(\dfrac{u^2_x(t,x)}2\Big)_x\,dx\,dt=
\dfrac{1}2\int_0^Ta(1)u^2_x(t,1)dt
-\dfrac{1}2\int_0^T\!\!\!\int_0^1 \big( xa(x)\big)' u^2_x(t,x)dx\,dt\,.
\end{equation}
Then the identity \eqref{eq:le1} follows by inserting \eqref{eq:l1_m1} and \eqref{eq:l1_e2} into \eqref{eq:le1_id1}.

Next, recall \eqref{eq:a} to obtain
\begin{eqnarray}\label{eq:anticipata}
\Big|\int_0^1 xu_x(t,x) u_t(t,x) dx\Big|\leqslant \frac12\int_0^1 \big\{u^2_t(t,x)+x^{2}u^2_x(t,x)\big\}  dx
\leqslant \dfrac{E_u(0)}{\min\{1, a(1)\}}\qquad\forall t\geqslant 0\,.
\end{eqnarray}
 Now, we deduce \eqref{eq:le0} from \eqref{eq:le1}, \eqref{eq:anticipata},  the inequality
$x|a'(x)|\leqslant 2 a(x)$, and the constancy of the energy. The conclusion has thus been proved for classical solutions.

In order  to extend  \eqref{eq:le0} and \eqref{eq:le1} to the mild solution associated with the initial data $(u_0,u_1)\in H^1_a(0,1)\times L^2(0,1)$, it suffices to approximate such data by $(u_0^n,u_1^n)\in H^2_a(0,1)\times H^1_a(0,1)$ and use 
\eqref{eq:le0} to show that the normal derivatives of the corresponding classical solutions give a Cauchy sequence in $L^2(0,T)$.
\end{Proof}

\begin{lemma}\label{le:2}
For any mild solution  $u$ of \eqref{DW} we have that, for every $T\geqslant 0$,
 \begin{equation}\label{eq:le2}
 \int_0^T\!\!\!\int_0^1\big\{a(x)u_x^2(t,x)-u_t^2(t,x)\big\} dtdx  + \Big[\int_0^1u(t,x) u_t(t,x) dx\Big]_{t=0}^{t=T}=0\,.
\end{equation}
\end{lemma}
\begin{Proof}
Once again we suppose $u$ is a classical solution of \eqref{DW}. Multiplying  equation \eqref{DW} by $u$ and integrating over $]0,T[\times]0,1[$ we obtain
\begin{eqnarray}
\nonumber
0&=&\int_0^T\!\!\!\int_0^1 u(t,x)\Big(u_{tt}(t,x) - \big(a(x)u_x(t,x)\big)_x\Big)dx\,dt
\\
\nonumber
&=&\Big[\int_0^1 u(t,x) u_t(t,x) dx\Big]_{t=0}^{t=T}-\int_0^T\!\!\!\int_0^1u^2_{t}(t,x)\,dx\,dt
\\
\nonumber
& &\qquad-\int_0^T\Big[a(x)\,u(t,x)u_{x}(t,x)\Big]_{x=0}^{x=1}dx\,dt
+\int_0^T\!\!\!\int_0^1a(x)\,u_x^2(t,x)dx\,dt\,.
\end{eqnarray}
The conclusion follows from the above identity because $a(x)\,u(t,x)u_{x}(t,x)$ vanishes at $x=1$ and, owing to \eqref{eq:H2a_2}, also at $x=0$. An approximation argument allows to extend the conclusion to mild solutions. \end{Proof}

\begin{theorem}
 \label{th:obser}
Assume \eqref{eq:hp_a} and let $u$ be the mild solution of \eqref{DW}-\eqref{b+i}. Then, for every $T\geqslant 0$,
\begin{equation}\label{eq:obser}
 a(1)\int_0^T u_x^2(t,1)\,dt\geqslant \Big\{(2-\mu_a)T-\dfrac{4}{\min\{1, a(1)\}}-2\,\mu_a\,\sqrt{C_a}\Big\}E_u(0)\,,
\end{equation}
where $C_a$ is the constant in \eqref{eq:cost_pcr}. 
\end{theorem}
\begin{Proof}
Suppose $u$ is a classical solution of \eqref{DW} (the general case can as usual be recovered  by an approximation argument).
By adding to the right-hand side of \eqref{eq:le1} the left side of \eqref{eq:le2} multiplied by $\mu_a/2$, we obtain
\begin{eqnarray*}
a(1)\int_0^T u_x^2(t,1)\,dt&=& \int_0^T\!\!\!\int_0^1\Big\{\Big(1-\dfrac{\mu_a}2\Big)u_t^2(t,x) + \Big[\Big(1+\dfrac{\mu_a}2\Big)a(x)-xa'(x)\Big]u_x^2(t,x)\Big\} dtdx
\\
& &\quad+ 2\Big[\int_0^1 xu_x(t,x) u_t(t,x) dx\Big]_{t=0}^{t=T}
+ \dfrac{\mu_a}2\Big[\int_0^1u(t,x) u_t(t,x) dx\Big]_{t=0}^{t=T}
\\
&\geqslant& (2-\mu_a)TE_u(0)+ 2\Big[\int_0^1 xu_x(t,x) u_t(t,x) dx\Big]_{t=0}^{t=T}
+ \dfrac{\mu_a}2\Big[\int_0^1u(t,x) u_t(t,x) dx\Big]_{t=0}^{t=T}\,,
\end{eqnarray*}
where we have use the inequality
$xa'(x)\leqslant \mu_a a(x)$ and the constancy of the energy. 
The conclusion follows from the above inequality recalling \eqref{eq:anticipata} and observing that 
\begin{eqnarray*}
\dfrac{1}2\Big|\int_0^1u(t,x) u_t(t,x) dx\Big|\leqslant 
\dfrac{1}2\int_0^1
\Big(\dfrac1{\sqrt{C_a}} \,u^2(t,x)+\sqrt{C_a}\,u_t^2(t,x)\Big)dx
\leqslant
\sqrt{C_a}\,E_u(0),
\end{eqnarray*}
where $C_a$ is  Poincar\'e's constant in \eqref{eq:cost_pcr}.
\end{Proof}

We recall that \eqref{DW} is said to be {\em observable} (via the normal derivative at $x=1$) {\em in time} $T>0$ if there exists a constant $C>0$ such that for any $(u_0,u_1)\in H^1_a(0,1)\times L^2(0,1)$  the mild solution of \eqref{DW}-\eqref{b+i} satisfies
\begin{equation}
\label{eq:defob}
 \int_0^T u_x^2(t,1)\,dt\geqslant C\,E_u(0)\,.
\end{equation}
Any constant satisfying \eqref{eq:defob} is called an {\em observability constant} for \eqref{DW} in time $T$. The supremum of all  observability constants for \eqref{DW}  is denoted by $C_T$. Equivalently, \eqref{DW} is observable if
\begin{equation*}
%\label{eq:}
C_T=\inf_{(u_0,u_1)\neq(0,0)}\dfrac{\int_0^T u_x^2(t,1)\,dt}{E_u(0)}>0\,.
\end{equation*}
The inverse $c_T=1/C_T$ is sometimes called the  cost of observability (or the  cost of control) in time $T$.
\begin{corollary}
Assume \eqref{eq:hp_a}. Then \eqref{DW} is observable  in time $T$ provided that
\begin{equation*}
%\label{eq:}
T>T_a:=\dfrac4{(2-\mu_a)\min\{1, a(1)\}}+2\,\mu_a\,\sqrt{C_a}\,,
\end{equation*} 
where $C_a$ is defined in \eqref{eq:cost_pcr}. In this case
\begin{equation*}
%\label{eq:}
C_T\geqslant \dfrac1{a(1)}\Big\{(2-\mu_a)T-\dfrac{4}{\min\{1, a(1)\}}-2\,\mu_a\,\sqrt{C_a}\Big\}.
\end{equation*}
\end{corollary}

\begin{remark}\label{re:observ_power}\rm Let $a$ be any of the two functions in Example~\ref{exa:power}. Then we can apply the above to conclude that, defining 
\begin{equation}
\label{eq:Ttheta}
T_\theta=\dfrac1{2-\theta}\Big(4+2\theta\,\min \Big\{ 2, \dfrac{1}{ \sqrt{2-\theta}}\Big\}\Big),
\end{equation}
we have that
 \begin{equation}\label{eq:observ_power}
C_T\geqslant (2-\theta)(T-T_\theta)\qquad \forall\,T\geqslant T_\theta\,.
\end{equation}
Observe that  $T_\theta\to 2$ as $\theta\downarrow 0$, which coincides with the classical observability time for the wave equation.
\end{remark}

\subsection{Failure of boundary observability}\label{se:failure}
In this section, we shall see that boundary observability is no longer  true when the constant $\mu_a$ in \eqref{eq:hp_a} is greater than or equal to $2$ and that, for $\mu_a<2$, the controllability time blows up as $\mu_a\uparrow 2$. We discuss two examples with power-like coefficients.
\begin{example}\rm\label{alpha=2}
Given $T>0$, consider the problem
 \begin{equation}
\label{DW2}
\begin{cases}
u_{tt} - \big(x^2u_x\big)_x=0  &\mbox{in}\;\; ]0,T[\times]0,1[
\\
\mbox{boundary conditions:}\;\;
u(t,1)=0 \;\;\mbox{and}\;\;
\lim_{x\downarrow 0} 
x^2\,u_x(t,x)=0
&0<t<T
\\
\mbox{initial conditions:}\;\;
\begin{cases}
\hspace{.cm}u(0,x)=u_0(x)
\\
u_t(0,x)= u_1(x)
\hspace{.cm}
\end{cases}
& x\in]0,1[\,,
\end{cases}
\end{equation}
where $u_0$ and $u_1$ are smooth functions with compact support in $]0,1[$, not identically zero. Observe that the  so-called Liouville transform
\begin{equation*}
%\label{eq:}
u(t,x)=\dfrac 1{\sqrt{x}}v\Big(t,\log \frac 1x\Big)
\end{equation*}
turns problem \eqref{DW2} into
 \begin{equation}
\label{W2}
\begin{cases}
v_{tt} - v_{yy}+\dfrac 14 v=0  &\mbox{in}\;\; ]0,T[\times]0,\infty[
\\
v(t,0)=0
&0<t<T
\\
\mbox{initial conditions:}\;\;
\begin{cases}
\hspace{.cm}v(0,y)=e^{-y/2}u_0(e^{-y}):=v_0(y)
\\
v_t(0,y)= e^{-y/2}u_1(e^{-y}):=v_1(y)\,.
\hspace{.cm}
\end{cases}
& y\in]0,\infty[\,.
\end{cases}
\end{equation}
Notice that $v_0$ and $v_1$ are, in turn, smooth functions with compact support in $]0,\infty[$ and the extreme point $y=0$ of the $y$-domain corresponds to $x=1$ in the $x$-domain. Since, for the wave equation (with a bounded potential) the support of the initial data propagates at finite speed (see, for instance, \cite{Evans}), the normal derivative $v_y(\cdot,0)$ of the solution to \eqref{W2} may well be identically zero on $[0,T]$ when the support of $v_0$ and $v_1$ is sufficiently far from $y=0$. Consequently, problem \eqref{DW2} is not observable on $[0,T]$ via the normal derivative $u_x(\cdot,1)$.
% where $\varphi:]0,1]\to[0,\infty[$ is defined by
%\begin{equation*}
%%\label{eq:}
%\varphi(x)=\int_x^1\dfrac{ds}{s}=\,,
%\end{equation*}

\end{example}

\begin{example}\rm\label{alpha>2}
Given $T>0$ and $\theta>2$, consider the problem
 \begin{equation}
\label{DWtheta}
\begin{cases}
u_{tt} - \big(x^\theta u_x\big)_x=0  &\mbox{in}\;\; ]0,T[\times]0,1[
\\
\mbox{boundary conditions:}\;\;
u(t,1)=0 \;\;\mbox{and}\;\;
\lim_{x\downarrow 0} 
x^\theta\,u_x(t,x)=0
&0<t<T
\\
\mbox{initial conditions:}\;\;
\begin{cases}
\hspace{.cm}u(0,x)=u_0(x)
\\
u_t(0,x)= u_1(x)
\hspace{.cm}
\end{cases}
& x\in]0,1[\,,
\end{cases}
\end{equation}
where $u_0$ and $u_1$ are smooth functions with compact support in $]0,1[$. Define $\varphi:]0,1]\to[0,\infty[$ by
\begin{equation*}
%\label{eq:}
\varphi(x)=\int_x^1\dfrac{ds}{s^{\theta/2}}=\dfrac{2(x^{1-\theta/2}-1)}{\theta-2}\,,
\end{equation*}
and denote by $\psi$ the inverse of $\varphi$, that is,
\begin{equation*}
%\label{eq:}
\psi(y)=\Big(\dfrac{2}{2+(\theta-2)y}\Big)^{\frac2{\theta-2}}\qquad\forall y\in [0,\infty[\,.
\end{equation*}
As in Example~\ref{alpha=2}, the change of variable 
\begin{equation*}
%\label{eq:}
u(t,x)=\dfrac 1{x^{\theta/4}}v\big(t,\varphi(x)\big)
\end{equation*}
transforms problem \eqref{DWtheta} into
 \begin{equation}
\label{Wtheta}
\begin{cases}
v_{tt} - v_{yy}+\,\dfrac{c(\theta)}{[2+(\theta-2)y]^2}\,v =0  &\mbox{in}\;\; ]0,T[\times]0,\infty[
\\
v(t,0)=0
&0<t<T
\\
\mbox{initial conditions:}\;\;
\begin{cases}
\hspace{.cm}v(0,y)=v_0(y)
\\
v_t(0,y)=v_1(y)
\hspace{.cm}
\end{cases}
& y\in]0,\infty[\,,
\end{cases}
\end{equation}
where $c(\theta)=\theta(3\theta-4)/4$, 
\begin{equation*}
%\label{eq:}
v_0(y)=\psi(y)^{\theta/4}u_0\big(\psi(y)\big)\,,\quad\mbox{and}\quad v_1(y)=\psi(y)^{\theta/4}u_1\big(\psi(y)\big)\,.
\end{equation*}
Notice that, as before, $v_0$ and $v_1$ are smooth functions with compact support in $]0,\infty[$ and the extreme point $y=0$ of the $y$-domain corresponds to $x=1$ in the $x$-domain. Therefore, the finite speed of propagation of the support for the wave equation (with a bounded potential) implies that the normal derivative $v_y(\cdot,0)$  is identically zero on $[0,T]$ when the support of $v_0$ and $v_1$ is sufficiently far from $y=0$. Consequently, problem \eqref{DW2} is not observable on $[0,T]$.
\end{example}

\subsection{Blow-up of observability time}
In this section, we will show that, for any fixed $T>0$ the observability constant $C_T(\theta)$ of \eqref{DWtheta},
with $0\leqslant \theta<2$,
goes to zero as $\theta\uparrow 2$. We begin by recalling spectral results for the family of Sturm-Liouville eigenvalue problems
\begin{equation}
\label{SL}
\begin{cases}
- \big(x^\theta y'(x)\big)'=\lambda y(x)  &x\in]0,1[
\\
\lim_{x\downarrow 0} 
x^\theta\,y'(x)=0 \;\;\mbox{and}\;\;
y(1)=0\,.
\end{cases}
\end{equation}
For any $\nu\geqslant 0$, denote by $J_\nu$ the Bessel function of the first kind of order $\nu$, that is,
\begin{equation*}
%\label{eq:}
J_\nu(x)=\sum_{m=0}^\infty\dfrac{(-1)^m}{m!\,\Gamma(m+\nu+1)}\Big(\dfrac x2\Big)^{2m+\nu}\qquad(x\geqslant 0),
\end{equation*}
where $\Gamma$ is Euler's Gamma function. Let $j_\nu$ be the first positive zero of $J_\nu$. 
\begin{proposition}
 \label{pr:SL}
 Given $\theta\in [1,2[$ define
 \begin{equation}
\label{eq:nu_kappa}
\nu_\theta=\dfrac{\theta-1}{2-\theta}\quad\mbox{and}\quad \kappa_\theta=\dfrac{2-\theta}{2}.
\end{equation}
Then the first eigenvalue of \eqref{SL} is given by  $\lambda_\theta=\kappa_\theta^2j^2_{\nu_\theta}$
and the corresponding normalized eigenfunction is
\begin{equation*}
%\label{eq:}
y_\theta(x)=\dfrac{\sqrt{2\kappa_\theta}}{\big|J'_{\nu_\theta}(j_{\nu_\theta})\big|}\,x^{\frac{1-\theta}2}\,J_{\nu_\theta}\big(j_{\nu_\theta}x^{\kappa_\theta}\big)\qquad(0<x<1)\,.
\end{equation*}
\end{proposition}
See \cite{Watson} for the proof.
\begin{theorem}\label{th:upper_bound}
For any fixed $T>0$ the observability constant $C_T(\theta)$ of \eqref{DWtheta}, with $1\leqslant \theta<2$, satisfies
\begin{equation}\label{eq:upper_bound}
C_T(\theta)\leqslant (2-\theta)T\,.
\end{equation}
\end{theorem}
\begin{Proof} Define 
 \begin{equation*}
%\label{eq:}
u_\theta(t,x)=\sin\big(\sqrt{\lambda_\theta} t\big)\,y_\theta(x)\qquad (t,x)\in ]0,T[\times]0,1[.
\end{equation*}
Then $u_\theta$ satisfies \eqref{DWtheta} with $u_0\equiv 0$ and $u_1(x)=\sqrt{\lambda_\theta}\,y_\theta(x)$. Now, straightforward computations lead to
\begin{equation*}
%\label{eq:}
\dfrac{\int_0^T |\partial_xu_\theta|^2(t,1)\,dt}{E_{u_\theta}(0)}=2T\kappa_\theta\Big(1-\dfrac{\sin\big(2\sqrt{\lambda_\theta}T)}{2\sqrt{\lambda_\theta}T}\Big)<(2-\theta)T
\end{equation*}
taking into account the definition of $\kappa_\theta$ in \eqref{eq:nu_kappa}. The conclusion follows recalling the definition of $C_T$.
 \end{Proof}

\begin{remark}\rm
 Given $C>0$, let $T^*_\theta(C)$ denote the infimum of all times $T>0$ such that $C$ is an observability constant for \eqref{DWtheta} in time $T$.
 Then \eqref{eq:upper_bound} yields
 \begin{eqnarray*}
T^*_a(C)\geqslant \dfrac C{2-\theta}\,,
\end{eqnarray*}
which means that the observability time $T^*_\theta(C)$ blows up, as $\theta\uparrow 2$, at essentially the same speed as $T_\theta$ in \eqref{eq:Ttheta}.
\end{remark}
\subsection{Controllability}
We consider the following controlled degenerate system
\begin{equation}
\label{DWC}
y_{tt} - \big(a(x)y_x\big)_x=0  \quad \mbox{in}\;\; ]0,\infty[\times]0,1[
\end{equation}
with
\begin{equation}
\label{b+iC}
\begin{cases}
\mbox{boundary conditions}\;\;
y(t,1)=f \;\;\mbox{and}\;\;\begin{cases}
\hspace{.cm}
y(t,0)=0
&
\mbox{if}\;\;\mu_a\in[0,1[
\\
\hspace{.cm}
\lim_{x\downarrow 0} 
a(x)\,y_x(t,x)=0
&
\mbox{if}\;\;\mu_a\in[1,2[
\end{cases}
&0<t<\infty
\\
\mbox{initial conditions}\;\;
\begin{cases}
\hspace{.cm}y(0,x)=y_0(x)
\\
y_t(0,x)= y_1(x)
\hspace{.cm}
\end{cases}
& x\in]0,1[.
\end{cases}
\end{equation}
where $f \in L^2(0,T)$ is the control. The solution of this controlled system is defined by transposition. At this stage, we have to introduce some notation.
Let us define the operator $A_0: D(A_0) \subset H \mapsto H$ where $D(A_0)=H^2_a(0,1)$ and $A_0u=:-(au')'$ for $u \in D(A_0)$.
We define $H^{-1}_a(0,1)$ as the dual space of $H^1_a(0,1)$ with respect to the pivot space $L^2(0,1)$. Then, thanks to
Proposition~\ref{pr:poincare}, one can prove that $A_0$ is an isomorphism
from $H^1_a(0,1)$ onto $H^{-1}_a(0,1)$. In particular, we have $H^{-1}_a(0,1)=A_0H^1_a(0,1)$.

\begin{definition}\label{transpose}
Let $f \in L^2_{loc}(0,\infty)$ and let $(y_0,y_1) \in L^2(0,1) \times H^{-1}_a(0,1)$ be fixed arbitrarily. We say that $y$ is a solution by transposition
of \eqref{DWC}-\eqref{b+iC} if 
\begin{equation*}
%\label{eq:}
y\in   \mathcal C^1\big([0,\infty[;H^{-1}_a(0,1)\big)\cap \mathcal C\big([0,\infty[;L^2(0,1)\big)
\end{equation*}
%$(y,y') \in \mathcal{C}([0,\infty); L^2(0,1)\times H^{-1}_a(0,1))$ 
satisfies for all $T>0$
\begin{multline}\label{transp1}
\langle y'(T), w_T^0 \rangle_{ H^{-1}_a(0,1), H^1_a(0,1)}- \int_0^1y(T)w^1_T dx=
\langle y_1, w(0) \rangle_{ H^{-1}_a(0,1), H^1_a(0,1)}- \int_0^1y_0w'(0) dx \\
+ \int_0^T f(t) w_x(t,1)dt \quad \forall \  (w_T^0, w_T^1) \in H^1_a(0,1) \times L^2(0,1),
\end{multline}
where $w$ is the solution of the backward equation
\begin{equation}
\label{DWB}
w_{tt} - \big(a(x)w_x\big)_x=0  \quad \mbox{in}\;\; ]0,\infty[\times]0,1[
\end{equation}
with
\begin{equation}
\label{b+iB}
\begin{cases}
\mbox{boundary conditions}\;\;
w(t,1)=0 \;\;\mbox{and}\;\;\begin{cases}
\hspace{.cm}
w(t,0)=0
&
\mbox{if}\;\;\mu_a\in[0,1[
\\
\hspace{.cm}
\lim_{x\downarrow 0} 
a(x)\,w_x(t,x)=0
&
\mbox{if}\;\;\mu_a\in[1,2[
\end{cases}
&0<t<\infty
\\
\mbox{final conditions}\;\;
\begin{cases}
\hspace{.cm}w(T,x)=w_T^0(x)
\\
w_t(T,x)= w_T^1(x)
\hspace{.cm}
\end{cases}
& x\in]0,1[.
\end{cases}
\end{equation}
\end{definition}
Note that thanks to the change of variable $u(t,x)=w(T-t,x)$ and to our previous results, the backward problem \eqref{DWB} admits a unique solution
$w\in   \mathcal C^1\big([0,\infty[;L^2(0,1)\big)\cap \mathcal  C\big([0,\infty[;H^1_a(0,1)\big)$. Moreover, this solution depends continuously on $W^T=:(w_T^0, w_T^1) \in
H^1_a(0,1)\times L^2(0,1)$
and the energy $E_w$ of $w$ is conserved through time. Now thanks to the direct inequality \eqref{eq:le0}, we have
$$
\int_0^Tw_x^2(t,1)dt \leqslant D_T E_w(0)=D_T E_w(T).
$$
Thus, the right hand side of \eqref{transp1} defines a continuous linear form with respect to $(w_T^0, w_T^1) \in
H^1_a(0,1)\times L^2(0,1)$. Moreover, this linear form depends continuously on $T>0$, for all $T>0$. Therefore, there is a unique solution
by transposition $y \in  \mathcal C^1\big([0,\infty[;H^{-1}_a(0,1)\big)\cap \mathcal C\big([0,\infty[;L^2(0,1)\big)$ of \eqref{transp1}.

Let $(y_0,y_1) \in L^2(0,1)\times H^{-1}_a(0,1)$, $(y_0^T,y_1^T) \in L^2(0,1)\times H^{-1}_a(0,1)$ be given: then one wants to determine if there exists a control $f \in L^2(0,T)$ such that
the solution of \eqref{DWC} satisfies $(y,y_t)(T,\cdot)\equiv (y_0^T,y_1^T)(\cdot)$. If this is possible for every $(y_0,y_1) \in L^2(0,1)\times H^{-1}_a(0,1)$
and $(y_0^T,y_1^T) \in L^2(0,1)\times H^{-1}_a(0,1)$, one says that \eqref{DWC} is exactly controllable in $L^2(0,1)\times H^{-1}_a(0,1)$.

By linearity and reversibility, it is easy to check that this
property will hold as soon as it holds for arbitrary initial data $(y_0,y_1)$ and for a zero final state, that is for $(y_0^T,y_1^T)(\cdot)=(0,0)$.

Let us consider the bilinear form $\Lambda$ defined on $H^1_a(0,1) \times L^2(0,1)$ by
$$
\Lambda(W^T, \widetilde{W}^T)=:\int_0^T w_x(t,1)\widetilde{w}_x(t,1)dt \quad \forall \ W^T, \widetilde{W}^T \in H^1_a(0,1) \times L^2(0,1).
$$
Thanks to the direct inequality $\Lambda$ is continuous on $H^1_a(0,1) \times L^2(0,1)$. Moreover thanks to the observability inequality \eqref{eq:obser}
$\Lambda$ is coercive on $H^1_a(0,1) \times L^2(0,1)$ for $T>T_a$. We also define the continuous linear map
$$
\mathcal{L}(W^T):= \langle y_1, w(0) \rangle_{ H^{-1}_a(0,1), H^1_a(0,1)}- \int_0^1 y_0 w'(0)dx  \quad \forall \ W^T \in H^1_a(0,1) \times L^2(0,1).
$$
Since $\Lambda$ is continuous and coercive on $H^1_a(0,1) \times L^2(0,1)$, and $\mathcal{L}$ is continuous on the Hilbert space
$H^1_a(0,1) \times L^2(0,1)$, we can apply the Lax-Milgram Lemma. This implies that there exists a unique $W^T \in H^1_a(0,1) \times L^2(0,1)$ such that
$$
\Lambda(W^T, \widetilde{W}^T)=-\mathcal{L}(\widetilde{W}^T). 
$$
We set $f=w_x(t,1)$ and denote by $y$ the solution by transposition of \eqref{DWC}. Then we have
\begin{multline*}
\int_0^T f(t)\tilde{w}_x(t,1)dt=\int w_x(t,1) \tilde{w}_x(t,1)dt=\Lambda(W^T, \widetilde{W}^T)= -\langle y_1, \widetilde{w}(0) \rangle_{ H^{-1}_a(0,1), H^1_a(0,1)}+\int_0^1 y_0 \widetilde{w}'(0)dx\\
\quad \forall \  (\widetilde{w}_T^0, \widetilde{w}_T^1) \in H^1_a(0,1) \times L^2(0,1).
\end{multline*}
On the other hand, by definition of the transposition solutions, we have
\begin{multline}\label{transp2}
 \int_0^T f(t) \widetilde{w}_x(t,1)dt=
\langle y'(T), \widetilde{w}_T^0 \rangle_{ H^{-1}_a(0,1), H^1_a(0,1)}- \int_0^1y(T)\widetilde{w}^1_T dx-
\langle y_1, \widetilde{w}(0) \rangle_{ H^{-1}_a(0,1), H^1_a(0,1)}+ \int_0^1y_0\widetilde{w}'(0) dx +\\
 \quad \forall \  (\widetilde{w}_T^0, \widetilde{w}_T^1) \in H^1_a(0,1) \times L^2(0,1),
\end{multline}
Hence, comparing these two last relations, we deduce that
$$
\langle y'(T), \widetilde{w}_T^0 \rangle_{ H^{-1}_a(0,1), H^1_a(0,1)}- \int_0^1y(T)\widetilde{w}^1_T dx=0 \quad \forall \  (\widetilde{w}_T^0, \widetilde{w}_T^1) \in H^1_a(0,1) \times L^2(0,1).
$$
Thus, we have
$$
(y,y')(T,\cdot)\equiv (0,0) \mbox{ on  }(0,1).
$$
%%%%%%%%%%%%%%%%%%%%%%%%%%%%%%%%%%%%%%%%%%%%%%%%%%%%%%%%%%%%%%%%
\section{Stabilization}
%%%%%%%%%%%%%%%%%%%%%%%%%%%%%%%%%%%%%%%%%%%%%%%%%%%%%%%%%%%%%%%%

\subsection{Linear stabilization}
Given $a$ satisfying assumptions \eqref{eq:hp_a}, let $\mu_a\in [0,2[$ be the constant in assumption $(ii)$. Consider the degenerate wave equation with boundary damping
\begin{equation}
\label{DWs}
u_{tt} - \big(a(x)u_x\big)_x=0  \quad \mbox{in}\;\; ]0,T[\times]0,1[
\end{equation}
with
\begin{equation}
\label{b+is}
\begin{cases}
\;\;
u_t(t,1)+u_x(t,1)+\beta u(t,1)=0 \,,\quad
%\\
%\mbox{and}\;\;
\begin{cases}
\hspace{.cm}
u(t,0)=0
&
\mbox{if}\;\;\mu_a\in[0,1[
\\
\hspace{.cm}
\lim_{x\downarrow 0} 
a(x)\,u_x(t,x)=0
&
\mbox{if}\;\;\mu_a\in[1,2[
\end{cases}
&(0<t<T)
\\
u(0,x)=u_0(x)\,,\quad
u_t(0,x)= u_1(x)
& (0\leqslant x\leqslant 1).
\end{cases}
\end{equation}
where $\beta\geqslant 0$ is given.

\subsection{Well-posedness}
Let us denote by $W^1_a(0,1)$ the space $V^1_a(0,1)$ itself, if $\mu_a\in [1,2[$, and the closed subspace of $V^1_a(0,1)$ consisting of all the functions $u\in V^1_a(0,1)$ such that $u(0)=0$,  if $\mu_a\in [0,1[$. Moreover, we set 
$$W^2_a(0,1)= V^2_a(0,1)\cap W^1_a(0,1).$$
Notice that $W^2_a(0,1)=V^2_a(0,1)$ when $\mu_a\in [1,2[$.

Now, consider the Hilbert space $\Hsta=W^1_a(0,1)\times L^2(0,1)$ with the scalar product
\begin{equation*}
%\label{eq:}
\big\langle (u,v),(\widetilde u, \widetilde v)\big\rangle=
\int_0^1\big(v(x)\widetilde v(x)+a(x)u'(x)\widetilde u'(x)\big)dx+a(1)\beta u(1)\widetilde u(1)
\qquad\forall\,(u,v),(\widetilde u,\widetilde v)\in \Hsta
\end{equation*}
and the unbounded operator $A_\beta:D(A_\beta)\subset \Hsta\to \Hsta$ defined by
\begin{equation*}
%\label{eq:}
\begin{cases}
\hspace{.cm}
D(A_\beta)=\big\{(u,v)\in W^2_a(0,1)\times W^1_a(0,1)~:~u'(1)+v(1)+\beta u(1)=0 \big\}
&
\\
A_\beta(u,v)=\big(v, (au')'\big)\qquad\forall (u,v)\in D(A_\beta)\,.
\hspace{.cm}
\end{cases}
\end{equation*}
Observe that $u'(1),v(1),$ and $\beta u(1)$ are well defined for all $(u,v)\in W^2_a(0,1)\times W^1_a(0,1)$ because of the classical Sobolev embedding theorem.
\begin{proposition}\label{pr:swell-posed}
Assume \eqref{eq:hp_a}. Then $A_\beta$ is a maximal dissipative operator on $\Hsta$.
\end{proposition}
\begin{Proof}
Let $(u,v)\in D(A_\beta)$. Then
\begin{eqnarray*}
\big\langle A_\beta(u,v),( u,  v)\big\rangle&=&\int_0^1\big((au')'v+au'v'\big)dx+a(1)\beta u(1)v(1)
\\
&=&a(1) v(1)\big(u'(1)+\beta u(1)\big)=-a(1) v^2(1)\leqslant 0\,.
\end{eqnarray*}
Therefore, $A_\beta$ is dissipative.

In order to show that $A_\beta$ is maximal dissipative, it remains to check that $I-A_\beta$ is onto. Equivalently,
given any $(f,g)\in \Hsta$, we have to solve the problem
\begin{equation}
\label{eq:max}
\begin{cases}
\hspace{.cm}
(u,v)\in D(A_\beta)
&
\\
v=u-f
\\
u-(au')'=f+g\,.
\hspace{.cm}
\end{cases}
\end{equation}
Consider the bilinear form $b:W^1_a(0,1)\times W^1_a(0,1)\to\R$ given by
\begin{equation*}
%\label{eq:}
b(u,\phi)=\int_0^1\big(u\phi+au'\phi'\big)dx +(\beta+1)a(1)u(1)\phi(1)\,,
\end{equation*}
and the linear form $L: W^1_a(0,1) \to\R$ given by
\begin{equation*}
%\label{eq:var}
L\phi=\int_0^1(f+g)\phi dx +a(1)\phi(1)f(1)\,.
\end{equation*}
In view of Proposition~\ref{pr:sobolev}, $b$ is a continuous bilinear form on $W^1_a(0,1)\times W^1_a(0,1)$ and $L$ is a continuous linear functional on $W^1_a(0,1)$. Moreover since $\beta \geqslant 0$, $b$ is also coercive on $W^1_a(0,1)\times W^1_a(0,1)$. So, by the Lax-Milgram Theorem there exists a unique solution $u \in W^1_a(0,1)$ of the variational problem
\begin{equation}
\label{eq:var}
b(u,\phi)=L\phi \qquad \forall \ \phi \in W^1_a(0,1)\,.
\end{equation}
We prove that $(u,v) \in D(A_\beta)$ and solves \eqref{eq:max} as follows. We denote by $\mathcal{C}^{\infty}_c(0,1)$ the space of functions which are in
 $\mathcal{C}^{\infty}(0,1)$ with compact support in $(0,1)$. Since $\mathcal{C}^{\infty}_c(0,1) \subset W^1_a(0,1)$, we have
 $$
 \int_0^1\big(u\phi+au'\phi'\big)dx =\int_0^1(f+g)\phi dx  \qquad \forall \ \phi \in \mathcal{C}^{\infty}_c(0,1) \,.
 $$
 Hence by duality, we have $u-(au')'=f+g$ in the sense of distributions. Thus $u \in  W^2_a(0,1)$ and
 $$
 u-(au')'=f+g \qquad \mbox{a.e in } (0,1)\,.
 $$
 Thus, we deduce after an integration by parts together with \eqref{eq:H2a_2} that 
 $$
 \int_0^1 u\phi dx + \int_0^1 a u'\phi'dx - a(1)u'(1)\phi(1)=  \int_0^1(f+g)\phi dx \qquad \forall \ \phi \in W^1_a(0,1)\,.
 $$
 This combined with \eqref{eq:var} yields
 $$
 a(1)\phi(1)\big(u'(1) +(\beta +1)u(1) -f(1)\big)=0 \qquad \forall \ \phi \in W^1_a(0,1)\,.
 $$
 Since $a(1)>0$ and the function $\phi$ defined by $\phi(x)=x$ for all $x \in (0,1)$ is in $W^1_a(0,1)$ we deduce
 that 
 $$
 u'(1) +(\beta +1)u(1) -f(1)=0 \,.
 $$
Setting $v=u-f$, we check that $(u,v) \in D(A_\beta)$ and solves \eqref{eq:max}.
\end{Proof}
Therefore, $A_\beta$ is the generator of a contraction semigroup in $\Hsta$, denoted by $e^{tA_\beta}$. For any $U_0:=(u_0,u_1)\in \Hsta$, $U(t):= e^{tA_\beta}U_0$  can be viewed as the  weak solution of the Cauchy problem
\begin{equation}\label{eq:Cauchy}
\begin{cases}
\hspace{.cm}
U'(t)=A_\beta U(t)
&t>0
\\
\hspace{.cm}
U(0)=U_0.
\end{cases}
\end{equation}
Moreover, the above solution is classical when $U_0\in D(A_\beta)$. We thus have the following result.
\begin{corollary}
Assume \eqref{eq:hp_a}. Then, for any $U_0=(u_0,u_1)\in D(A_\beta)$, problem~\eqref{eq:Cauchy} has a unique solution
\begin{equation*}
%\label{eq:}
U\in C^1\big([0,\infty);\Hsta\big)\cap C\big([0,\infty);D(A_\beta)\big)
\end{equation*}
given by $U(t)=e^{tA_\beta}U_0$. Moreover, setting $U(t)=\big(u(t),v(t)\big)$, we have that 
\begin{itemize}
\item $u$ is the unique solution of problem~\eqref{DWs}-\eqref{b+is} such that
\begin{equation*}
%\label{eq:}
u\in C^2\big([0,\infty);L^2(0,1)\big)\cap C^1\big([0,\infty);W^1_a(0,1)\big)\cap C\big([0,\infty);W^2_a(0,1)\big),
\end{equation*}
\item the energy of $u$ defined by
\begin{equation}\label{energy_u}
E_u(t)=:\dfrac{1}{2}\Big[\int_0^1\Big( u_t^2 + a u_x^2\Big)dx + \beta a(1) u^2(t,1)
\Big] 
\end{equation}
satisfies
\begin{equation}\label{dissip_u}
\dfrac{d E_u}{dt}(t)=-a(1) u_t^2(t,1) \leqslant 0 \qquad \forall \ t \geqslant 0.
\end{equation}
\end{itemize}
\end{corollary}

%%%%%%%%%%%%%%%%%%%%%%%%%%%%%%%%%
We shall need the following results in the sequel.
\begin{proposition}\label{pr:poincare-w}
 Assume \eqref{eq:hp_a}.
 Then 
\begin{equation}\label{eq:poincare-w}
\|u\|^2_{L^2(0,1)}\leqslant 2 |u(1)|^2 + C'_a\; |u|^2_{1,a}\qquad\forall\,u\in W^1_{a}(0,1),
\end{equation}
where
\begin{equation}\label{eq:cost_pcr}
C'_a=\frac{1}{a(1)} \min \Big\{ 4, \dfrac{2}{ 2-\mu_a}\Big\}\,.
\end{equation}
Moreover assume that $\beta >0$. Then, denoting by $|||\cdot|||_{1,a}$ the norm defined by
$$
|||u|||_{1,a}=\left(|u|_{1,a}^2 + \beta a(1) u^2(1)\right)^{1/2} \qquad u \in W^1_a(0,1)\,,
$$
we have
\begin{equation}\label{trip0}
|||u|||_{1,a}^2 \geqslant  \alpha_a ||u||^2_{L^2(0,1)}
\qquad\forall\,u\in W^1_{a}(0,1)\,.
\end{equation}
where
$$
\alpha_a=\min\left(\frac{1}{C_a'}, \frac{\beta a(1)}{2}\right)>0\,.
$$
Moreover we also have
\begin{equation}\label{trip1}
\dfrac{\alpha_a}{\alpha_a +1} \left(||u||_{1,a}^2 + \beta a(1) u^2(1)\right) \leqslant |||u|||_{1,a}^2 \leqslant \gamma_a ||u||_{1,a}^2 \qquad\forall\,u\in W^1_{a}(0,1)\,,
\end{equation}
where
$$
\gamma_a= \max\left(2\beta a(1), 1+ \frac{2\beta}{2-\mu_a}\right)\,.
$$
\end{proposition}
\begin{Proof} Let $u\in W^1_{a}(0,1)$. We follow the proof of Proposition~\ref{pr:poincare}. We give two different bounds for $\|u\|^2_{L^2(0,1)}$ in terms of $|u|^2_{1,a}$ and $u^2(1)$. The conclusion \eqref{eq:poincare-w} will follow by taking the minimum of the two corresponding constants.

First, for any $x\in ]0,1]$ we have that
\begin{equation*}
%\label{eq:}
|u(x)-u(1)|=\Big|\int_x^1u'(s)ds\Big|\leqslant |u|_{1,a}\, \Big\{\int_x^1\dfrac{ds}{a(s)}\Big\}^{\frac12}.
\end{equation*}
Therefore, proceeding as in the proof of \eqref{eq:fubini}, we have
\begin{eqnarray}\label{eq:fubini-w}
\int_0^1|u(x)-u(1)|^2dx\leqslant  |u|^2_{1,a}\, \int_0^1dx\int_x^1\dfrac{ds}{a(s)} =
|u|^2_{1,a}\, \int_0^1\dfrac{s}{a(s)} ds \leqslant \dfrac{\;  |u|^2_{1,a}}{a(1)(2-\mu_a)}.
\end{eqnarray}
Since
$$
\int_0^1|u(x)|^2dx\leqslant 2 |u(1)|^2 + 2 \int_0^1|u(x)-u(1)|^2dx\,,
$$
we deduce by \eqref{eq:fubini-w} the first bound we mentioned above, that is,
 \begin{equation}\label{eq:poincarew1}
\|u\|^2_{L^2(0,1)}\leqslant 2 |u(1)|^2 + \dfrac{\; 2 |u|^2_{1,a}}{a(1)(2-\mu_a)}\qquad\forall\,u\in W^1_{a}(0,1)\,.
\end{equation}

%Therefore, proceeding as in the proof of \eqref{eq:fubini}, we have
%\begin{eqnarray}\label{eq:fubini-w}
%\int_0^1|u(x)|^2dx\leqslant 2 |u(1)|^2 + |u|^2_{1,a}\, \int_0^1dx\int_x^1\dfrac{ds}{a(s)} =2 |u(1)|^2 +
%|u|^2_{1,a}\, \int_0^1\dfrac{s}{a(s)} ds\,.
%\end{eqnarray}
%By \eqref{eq:a} and \eqref{eq:fubini-w} we obtain the first bound we mentioned above, that is,
% \begin{equation}\label{eq:poincarew1}
%\|u\|^2_{L^2(0,1)}\leqslant 2 |u(1)|^2 + \dfrac{\; 2 |u|^2_{1,a}}{a(1)(2-\mu_a)}\qquad\forall\,u\in W^1_{a}(0,1).
%\end{equation}
%
Next, observe that, for all $x\in]0,1[$,
\begin{eqnarray*}
0&\leqslant&\int_x^1\Big(su'(s)+\frac12 u(s)\Big)^2ds
\\
&=&\int_x^1\Big(s^2|u'(s)|^2+\frac14 |u(s)|^2+s\,u(s)\,u'(s)\Big) ds
\\
&=&\int_x^1\Big(s^2|u'(s)|^2-\frac14 |u(s)|^2 \Big) ds+ \frac12 u^2(1)-\frac12 x\,|u(x)|^2.
\end{eqnarray*}
Therefore, taking the limit as $x\downarrow 0$, by \eqref{eq:a} and \eqref{eq:H1a_0} we obtain the announced second  bound:
\begin{equation}
\label{eq:poincarew2}
\int_0^1|u(s)|^2 ds\leqslant 2 u^2(1)+ 4\int_0^1 s^2|u'(s)|^2ds \leqslant 2 u^2(1)+\frac4{a(1)}\int_0^1 a(s)|u'(s)|^2ds\qquad\forall\,u\in W^1_{a}(0,1).
\end{equation}
The inequality \eqref{eq:poincare-w} follows from \eqref{eq:poincarew1} and \eqref{eq:poincarew2}.

We have
$$
|||u|||_{1,a}^2 \geqslant \min \left(\frac{1}{C_a'}, \frac{\beta a(1)}{2}\right) \left(2u^2(1) + C_a'|u|_{1,a}^2\right) \geqslant \alpha_a ||u||^2_{L^2(0,1)}
\qquad\forall\,u\in W^1_{a}(0,1)\,.
$$
Writing $1= \frac{\alpha_a}{\alpha_a+1} + \frac{1}{\alpha_a+1}$ and using the above inequality, we obtain
$$
|||u|||_{1,a}^2\geqslant  \frac{\alpha_a}{\alpha_a+1} |||u|||_{1,a}^2+ \frac{\alpha_a}{\alpha_a+1} ||u||^2_{L^2(0,1)} \qquad\forall\,u\in W^1_{a}(0,1)\,.
$$
This gives the left hand side of \eqref{trip1}.
On the other hand, since
$$
|u(1)|^2 \leqslant 2 \int_0^1|u(x)|^2dx  + 2 \int_0^1|u(x)-u(1)|^2dx \leqslant 2 ||u||^2_{L^2(0,1)} + \dfrac{\; 2 |u|^2_{1,a}}{a(1)(2-\mu_a)}\qquad\forall\,u\in W^1_{a}(0,1)\,.
$$
This inequality yields
$$
|||u|||_{1,a}^2 \leqslant 2 \beta a(1) ||u||^2_{L^2(0,1)} + \big(1+ \frac{2\beta}{2-\mu_a}\big) |u|_{1,a}^2 \leqslant  \max\left(2\beta a(1), 1+ \frac{2\beta}{2-\mu_a}\right)
||u||_{1,a}^2  \qquad\forall\,u\in W^1_{a}(0,1)\,.
$$
This gives the right inequality in \eqref{trip1}.
\end{Proof}
\begin{proposition}\label{pr:zerord}
 Assume \eqref{eq:hp_a} and that $\beta > 0$ is given.
 Then the variational problem
 \begin{equation}
\label{eq:var0}
 \int_0^1 az'\phi'dx + \beta a(1) z(1) \phi(1)= \lambda a(1) \phi(1)  \qquad \forall \ \phi \in W^1_a(0,1)\,.
\end{equation}
 admits a unique solution $z \in W^1_a(0,1)$ which satisfies the elliptic estimates
 \begin{equation}\label{ellip1}
|||z|||_{1,a}^2 \leqslant  \frac{a(1)}{\beta} \lambda^2 \,, \quad ||z||_{L^2(0,1)}^2 \leqslant \frac{a(1)}{\beta \alpha_a} \lambda^2 \,.
 \end{equation}
 Moreover $z \in W^2_a(0,1)$ and solves
 \begin{equation}
 \label{eq: var0s}
 \begin{cases}
 -(az')'=0 \,,\\
 z'(1) +\beta z(1)=\lambda \,.
 \end{cases}
 \end{equation}
  \end{proposition}
 \begin{Proof}
 We denote by $\tilde{b}$ the bilinear form on $W^1_a(0,1)$ defined by
 $$
 \tilde{b}(z,\phi)=: \int_0^1 az'\phi'dx + \beta a(1) z(1) \phi(1) \qquad z, \phi \in W^1_a(0,1)\,.
 $$
 Thanks to Proposition \ref{pr:poincare-w} (see \eqref{trip1}), $\tilde{b}$ is a symmetric continuous and coercive bilinear form on $W^1_a(0,1)$ and the
 linear form $\tilde{L}$ defined by $\tilde{L}\phi=:\lambda a(1) \phi(1)$ for $\phi \in W^1_a(0,1)$ is continuous. Hence thanks to the Lax-Milgram's Theorem,
 the above variational problem admits a unique solution $z \in W^1_a(0,1)$. Hence we have
 $$
|||z|||_{1,a}^2= \tilde{b}(z,z)=\lambda a(1) z(1) \leqslant  \frac{\sqrt{a(1)}}{\sqrt{\beta}} |\lambda| |||z|||_{1,a}\,.
 $$
 Hence we have
 $$
 |||z|||_{1,a}^2 \leqslant  \frac{a(1)}{\beta} \lambda^2 \,.
 $$
 This, together with \eqref{trip0} yields
 $$
 ||z||_{L^2(0,1)}^2 \leqslant \frac{a(1)}{\beta \alpha_a} \lambda^2 \,.
 $$
 Proceeding as in the proof of Proposition \ref{pr:swell-posed}, we
 show that $z \in W^2_a(0,1)$ and solves \eqref{eq: var0s}.
  \end{Proof}
  \begin{theorem}\label{Theorem_stab}
  Assume \eqref{eq:hp_a} and that $\beta > 0$ is given. Then for any $(u_0,u_1) \in \mathcal{H}_{\beta}$, the solution of \eqref{DWs}-\eqref{b+is}.
  satisfies the uniform exponential decay
  \begin{equation}\label{exp_decay}
  E_u(t)\leqslant E_u(0)e^{1-t/M_{a,\beta}},\quad\forall  t\in [M_{a,\beta},+\infty).
  \end{equation}
  where $M_{a,\beta}>0$ is given in \eqref{Mabeta} and is independent of $(u_0,u_1)$.
  \end{theorem}
  \begin{Proof}
 Let $U_0=(u_0,u_1)\in D(A_\beta)$ be given, and $U$ be the corresponding solution of problem~\eqref{eq:Cauchy}. Then we recall that setting as above $U(t)=\big(u(t),v(t)\big)$, we have that $u$ is the solution of problem~\eqref{DWs}-\eqref{b+is}. We multiply \eqref{DWs} by $xu_x$ and integrate the resulting equation over $(S,T) \times (0,1)$. This gives after suitable integrations by parts
 $$
 \int_S^T\int_0^1 \Big(-x \big(\dfrac{u_t^2}{2}\big)_x + a(x) u_x^2 + xa(x) \big(\dfrac{u_x^2}{2}\big)_x\Big)dxdt + \Big[\int_0^1xu_xu_t dx
 \Big]_S^T -\int_S^T \big[xa u_x^2\big]_0^1 dt=0 \quad \forall \ 0 \leqslant S \leqslant T.
 $$
 We integrate by parts twice again. This gives
 $$
 \int_S^T\int_0^1 \Big(\dfrac{u_t^2}{2} + (a-x a') \dfrac{u_x^2}{2}\Big)dxdt + \Big[\int_0^1xu_xu_t dx
 \Big]_S^T -\dfrac{1}{2} \int_S^T \Big(\big[xa u_x^2\big]_0^1 + \big[x u_t^2\big]_0^1\Big)dt=0 \quad \forall \ 0 \leqslant S \leqslant T.
 $$
 Now, we recall that $u(t,.) \in W^2_a(0,1)$ and $u_t(t,.) \in W^1_a(0,1)$ for every $t \geqslant 0$. Hence, since $W^p_a(0,1) \subset V^p_a(0,1)$ for
 $p=1,2$ and thanks to Proposition~\ref{pr:sobolev} (properties {\bf (I)} and {\bf (II)}), we have
 $$
 (xu_t^2(t,x))_{|_{x=0}}=0\,, (xa(x)u_x^2(t,x))_{|_{x=0}=0}. %\red{dire ou varie $t$}
 $$
 Using these two relations in the above equation, we obtain
 \begin{equation}\label{stab1}
 \int_S^T\int_0^1 \Big(\dfrac{u_t^2}{2} + (a-x a') \dfrac{u_x^2}{2}\Big)dxdt + \Big[\int_0^1xu_xu_t dx
 \Big]_S^T -\dfrac{1}{2} \int_S^T \Big(a(1) u_x^2(t,1) + u_t^2(t,1)\Big)dt=0 \quad \forall \ 0 \leqslant S \leqslant T.
 \end{equation}
 We multiply \eqref{DWs} by $u$ and integrate the resulting equation over $(S,T) \times (0,1)$. This gives after a suitable integration by parts.
 $$
 \int_S^T\int_0^1 \Big(-u_t^2 + au_x^2\Big) dxdt + \Big[\int_0^1u_t udx
 \Big]_S^T - \int_S^T \big[au_xu\big]_0^1dt=0 \quad \forall \ 0 \leqslant S \leqslant T.
 $$
Using now Proposition~\ref{pr:sobolev} (see {\bf (III)}), we have
$$
(a(x)u(t,x)u_x(t,x))_{|_{x=0}}=0, %\red{dire ou varie $t$}
$$
 so that
 \begin{equation}\label{stab2}
 \int_S^T\int_0^1 \Big(-u_t^2 + au_x^2\Big) dxdt + \Big[\int_0^1u_t udx
 \Big]_S^T - \int_S^T a(1)u_x(t,1)u(t,1)dt=0 \quad \forall \ 0 \leqslant S \leqslant T.
 \end{equation}
 We now combine \eqref{stab1} multiplied by 2 with \eqref{stab2} multiplied $\dfrac{\mu_a}{2}$. This gives
 \begin{multline}\label{stab3}
 \int_S^T\int_0^1 \Big[(2-\mu_a)\dfrac{u_t^2}{2} + \big[2(a-x a')+ a\mu_a\big] \dfrac{u_x^2}{2}\Big]dxdt + \dfrac{2-\mu_a}{2}\beta a(1)\int_S^T u^2(t,1)dt=\\
 -2 \Big[\int_0^1xu_xu_t dx \Big]_S^T - \dfrac{\mu_a}{2}\Big[\int_0^1u_tu dx \Big]_S^T+
 \int_S^T h(t)dt \quad \forall \ 0 \leqslant S \leqslant T,
 \end{multline}
 where the function $h$ is given by
 \begin{equation}\label{defh}
 h(t)=(1+a(1))u_t^2(t,1) + a(1) \beta (1+\beta-\mu_a)u^2(t,1) + \big(2\beta-\dfrac{\mu_a}{2}\big)a(1) u_t(t,1)u(t,1) \quad t \in (S,T).
 \end{equation}
 By definition of $\mu_a$, we have
 $$
 (2-\mu_a)a \leqslant 2(a-x a')+ a\mu_a.
 $$
 This, together with \eqref{stab3}, gives
 \begin{multline}\label{stab4}
(2-\mu_a) \int_S^TE_u(t)dt \leqslant 
 - \Big[\int_0^12xu_xu_t  + \dfrac{\mu_a}{2} u_tu dx \Big]_S^T +
 \int_S^T h(t)dt \quad \forall \ 0 \leqslant S \leqslant T.
 \end{multline}
On the other hand, we have
\begin{equation}\label{stab5}
h(t) \leqslant \eta_1 u_t^2(t,1) + \eta_2 a(1)u^2(t,1) \quad \forall \  t \in (S,T),
\end{equation}
where
$$
\eta_1= \big(1 + \dfrac{3}{2} a(1)\big) , \eta_2= \Big[ \beta(1+\beta -\mu_a) +\dfrac{1}{2}\big(2\beta -\dfrac{\mu_a}{2}\big)^2\Big].
$$
We also have
$$
\int_0^1 \Big|2xu_xu_t  + \dfrac{\mu_a}{2} u_tu\Big|dx \leqslant \int_0^1\Big[ x^2 u_x^2 + (1+ \dfrac{\mu_a}{4})u_t^2 + \dfrac{\mu_a}{4}u^2\Big]dx.
$$
Using \eqref{eq:late} together with \eqref{eq:poincare-w}, we deduce that
\begin{multline*}
\int_0^1 \Big|2xu_xu_t  + \dfrac{\mu_a}{2} u_tu\Big|dx \leqslant \int_0^1\Big[(1+ \dfrac{\mu_a}{4})u_t^2+ \big(\dfrac{1}{a(1)} +
\dfrac{\mu_a}{4}C_a'\big) au_x^2 + \dfrac{\mu_a}{2}u^2(1)\Big]dx \leqslant \
C_a''E_u(t) \quad \forall \  t \in [S,T],
\end{multline*}
where
$$
C_a''=2\max\Big(1+ \dfrac{\mu_a}{4}, \dfrac{1}{a(1)} +\dfrac{\mu_a}{4}C_a', \dfrac{\mu_a}{2\beta a(1)}
\Big).
$$
Using this inequality together with \eqref{stab5} in \eqref{stab4}, we obtain
\begin{equation}\label{stab6}
(2-\mu_a) \int_S^TE_u(t)dt \leqslant C_a''\Big(E_u(S)+E_u(T)\Big) + \eta_1\int_S^T u_t^2(t,1)dt + \eta_2 \int_S^T a(1) u^2(t,1)dt. 
\end{equation}
Using the dissipation relation \eqref{dissip_u}, we deduce that
\begin{multline}\label{stab7}
(2-\mu_a) \int_S^TE_u(t)dt \leqslant C_a''\Big(E_u(S)+E_u(T)\Big) + \dfrac{\eta_1}{a(1)}\Big(E_u(S)-E_u(T)\Big)+ \eta_2 \int_S^Ta(1) u^2(t,1)dt \leqslant \\
\Big(2 C_a'' + \dfrac{\eta_1}{a(1)}\Big) E_u(S) + \eta_2 \int_S^Ta(1) u^2(t,1)dt .
\end{multline}
We now estimate the last term of this inequality as follows. Set $\lambda=u(t,1)$ and denote by $z$ the solution of the
degenerate elliptic problem \eqref{eq: var0s}. We multiply \eqref{DWs} by $z$ and integrate the resulting equation over $(S,T) \times (0,1)$. This gives after suitable integrations by parts.
\begin{equation}\label{stab8}
\int_S^Ta(1)u^2(t,1)dt=\int_S^T\int_0^1 u_tz_t dxdt -a(1)\int_S^T u_t(t,1)z(t,1)dt-\Big[\int_0^1u_t zdx\Big]_S^T.
\end{equation}
We now estimate the terms of the right hand side in this inequality, as follows.
First, thanks to the second inequality in \eqref{ellip1}, we have
\begin{equation}\label{stab9}
||z_t||_{L^2(0,1)}^2 \leqslant \frac{a(1)}{\beta \alpha_a} u_t(t,1)^2.
\end{equation}
Moreover, thanks to the first inequality in \eqref{ellip1} and to the definition of $|||\cdot|||_{1,a}$, we have
$$
\beta a(1) z^2(t,1) \leqslant |||z|||_{1,a}^2 \leqslant  \frac{a(1)}{\beta} u^2(t,1),
$$
so that
\begin{equation}\label{stab10}
z^2(t,1) \leqslant  \frac{1}{\beta^2} u^2(t,1) \leqslant \dfrac{2}{\beta^3 a(1)} E_u(t).
\end{equation}
On the other hand, we have, thanks to the second inequality in \eqref{ellip1}
\begin{equation}\label{stab11}
\Big|\int_0^1u_t (t,x)z(t,x)dx \Big| \leqslant \dfrac{1}{\beta \sqrt{\alpha_a}}\left(\int_0^1 \dfrac{u_t^2}{2}dx + \dfrac{\beta a(1)}{2}  u^2(t,1)\right) \leqslant
 \dfrac{1}{\beta \sqrt{\alpha_a}}E_u(t) \quad \forall \ t \in [S,T].
\end{equation}
We now use \eqref{stab9}-\eqref{stab11} in \eqref{stab8}. This gives
$$
\int_S^Ta(1)u^2(t,1)dt \leqslant \delta \big(1 + \dfrac{1}{\beta^3}\big)\int_S^T E_u(t)dt + \dfrac{1}{2\delta}\Big(1+ \dfrac{1}{\beta \alpha_a}
\Big) \int_S^T a(1)u_t^2(t,1)dt +   \dfrac{1}{\beta \sqrt{\alpha_a}}\big(E_u(S)+ E_u(T)\big).
$$
Using now \eqref{dissip_u} in this estimate, we obtain
$$
\int_S^Ta(1)u^2(t,1)dt \leqslant \delta \big(1 + \dfrac{1}{\beta^3}\big)\int_S^T E_u(t)dt + \dfrac{1}{2\delta}\Big(1+ \dfrac{1}{\beta \alpha_a}
\Big) \big(E_u(S)- E_u(T)\big) +   \dfrac{1}{\beta \sqrt{\alpha_a}}\big(E_u(S)+ E_u(T)\big).
$$
We now choose $\delta =\dfrac{2-\mu_a}{2\eta_2 \big(1+\dfrac{1}{\beta^3}\big)}$ in the above inequality and combine the resulting
inequality in \eqref{stab7} to obtain
\begin{equation}\label{stab12}
\int_S^TE_u(t)dt \leqslant  M_{a, \beta} E_u(S),
\end{equation}
 where
 \begin{equation}\label{Mabeta}
 M_{a,\beta}= \dfrac{2}{(2-\mu_a)}\Big[2 C_a'' + \dfrac{\eta_1}{a(1)} + \dfrac{\eta_2^2\big(1+\dfrac{1}{\beta^3}\big)}{2-\mu_a}\Big(1+ \dfrac{1}{\beta \alpha_a}\Big) + \dfrac{2\eta_2}{\beta \sqrt{\alpha_a}}
 \Big].
 \end{equation}
 Now we use the following well-known result (see \cite[Theorem 8.1] {Komor94}).
 \begin{lemma}\label{lem-expo}
Assume that $ E: [0,+\infty) \mapsto [0,+\infty)$ is a non-increasing function and that  there is a constant  $M>0$ such that
   $$\int_t^\infty E(s)ds\leqslant M E(t), \quad\forall t\in  [0,+\infty).$$
Then we have
   $$E(t)\leqslant E(0)e^{1-t/M},\quad\forall  t\in [M,+\infty).$$
\end{lemma}
Applying this result on $E=E_u$ which is nonnegative, nonincreasing on $[0,\infty)$ and satisfies \eqref{stab12}, we have
$$
E_u(t)\leqslant E_u(0)e^{1-t/M_{a,\beta}},\quad\forall  t\in [M_{a,\beta},+\infty).
$$
  \end{Proof}
  
  \section{Nonlinear stabilization}
  
 In the previous section we considered the case of a linear boundary feedback. Here we extend our stability analysis to one-dimensional degenerate wave equations damped by a {\em nonlinear} boundary  feedback with arbitrary growth. For this, we combine our results for the linear case with the optimal-weight convexity method of \cite{alaamo2005, alajde2010}. 
  
  Let $\rho: \mathbb{R} \mapsto \mathbb{R}$ be a nondecreasing continuous function such that $\rho(0)=0$ and assume 
  there exist constants $c_1>0$, $c_2>0$ and an odd, continuously differentiable, strictly increasing function $g$ on $[-1,1]$ such that
 \begin{eqnarray}\label{hyprho}
c_1 g(|s|) \leqslant  | \rho(s)| \leqslant c_2g^{-1}(|s|) \quad \forall \  |s| \leqslant 1,\\
c_1 |s|  \leqslant  | \rho(s)| \leqslant c_2 |s| \quad \forall \  |s| \geqslant 1. \nonumber
  \end{eqnarray}
As before, let $a$ be given such that  assumptions \eqref{eq:hp_a} hold, and let $\mu_a\in [0,2[$ be the constant in assumption~$(ii)$. Consider the degenerate wave equation 
\begin{equation}
\label{DWsNL}
u_{tt} - \big(a(x)u_x\big)_x=0  \quad \mbox{in}\;\; ]0,T[\times]0,1[
\end{equation}
with the nonlinear boundary damping
\begin{equation}
\label{b+isNL}
\begin{cases}
\;\;
\rho(u_t(t,1))+u_x(t,1)+\beta u(t,1)=0 \,,\quad
%\\
%\mbox{and}\;\;
\begin{cases}
\hspace{.cm}
u(t,0)=0
&
\mbox{if}\;\;\mu_a\in[0,1[
\\
\hspace{.cm}
\lim_{x\downarrow 0} 
a(x)\,u_x(t,x)=0
&
\mbox{if}\;\;\mu_a\in[1,2[
\end{cases}
&(0<t<T)
\\
u(0,x)=u_0(x)\,,\quad
u_t(0,x)= u_1(x)
& (0\leqslant x\leqslant 1)
\end{cases}
\end{equation}
where $\beta\geqslant 0$ is given.

\begin{remark}\rm
Typical examples for $g$ are:
\begin{itemize}
\item  the linear case $g(x)=c x$ on $\mathbb{R}$,
\item  the polynomial case $g(x)=|x|^{p-1}x$  with $p>1$ in a neighborhood of $x=0$,
\item $g(x)= |x|^{p-1}x \ln^q(\frac{1}{|x|})$ with $p>1, q>0$ in a neighborhood of $x=0$, 
\item $g(x)=sign(x)e^{-1/x^2}$ in a neighborhood of $x=0$,
\item $g(x)=sign(x) e^{-\ln^p(\frac{1}{|x|})}$  with $1<p<2$ in a neighborhood of $x=0$.
\end{itemize}
See e.g. \cite{Komor94} for the linear and polynomial cases and \cite{alajde2010} for the other cases and the references therein, and \cite{alawanyu} for the last example when $p>2$.
\end{remark}

\subsection{Well-posedness}
We keep the functional spaces introduced in the previous section (for linear stabilization). However, we now need to deal with the nonlinear unbounded operator $A_\beta^{nl}:D(A_\beta^{nl})\subset \Hsta\to \Hsta$ defined by
\begin{equation*}
%\label{eq:}
\begin{cases}
\hspace{.cm}
D(A_\beta^{nl})=\big\{(u,v)\in W^2_a(0,1)\times W^1_a(0,1)~:~u'(1)+\rho(v(1))+\beta u(1)=0 \big\}
&
\\
A_\beta^{nl}(u,v)=\big(v, (au')'\big)\qquad\forall (u,v)\in D(A_\beta^{nl})\,.
\hspace{.cm}
\end{cases}
\end{equation*}
\begin{remark}\rm
Note that the set $D=\{(u,v)\in W^2_a(0,1)\times H^1_0(0,1)~:~u'(1)+\beta u(1)=0 \big\}$ is a subset of $D(A_\beta^{nl})$ and is dense
in $W^1_a(0,1) \times L^2(0,1)$. Therefore $D(A_\beta^{nl})$ is dense in $W^1_a(0,1) \times L^2(0,1)$.
\end{remark}

\begin{proposition}\label{pr:swell-posed}
Assume \eqref{eq:hp_a} and the above assumptions on $\rho$. Then $A_\beta^{nl}$ is a maximal dissipative operator on $\Hsta$.
\end{proposition}
\begin{Proof}
Let $(u,v), (w,z)\in D(A_\beta^{nl})$. Then
\begin{multline*}
\big\langle A_\beta^{nl}(u,v)- A_\beta^{nl}(w,z),( u,  v)-(w,z)\big\rangle=\int_0^1\big((a(u-w)')'(v-z)+a(u-w)'(v-z)'\big)dx+a(1)\beta (u-w)(1)(v-z)(1)=\\
a(1) (v-z)(1)\left[(u-w)'(1)+\beta (u-w)(1)\right]=-a(1) (v-z)(1) \left[\rho(v(1))-\rho(z(1))\right]\leqslant 0\,.
\end{multline*}
Therefore, $A_\beta^{nl}$ is dissipative. Let us now prove that $I-A_\beta^{nl}$ is onto. Equivalently,
given any $(f,g)\in \Hsta$, we have to solve the problem
\begin{equation}
\label{eq:maxNL}
\begin{cases}
\hspace{.cm}
(u,v)\in D(A_\beta^{nl})
&
\\
v=u-f
\\
u-(au')'=f+g\,.
\hspace{.cm}
\end{cases}
\end{equation}
Let us define
\begin{equation}\label{primrho}
R(s)=\int_0^s \rho(\tau)d\tau \quad \forall \ s \in \mathbb{R}.
\end{equation}
We define the functional $J_{\beta}: W^1_a(0,1) \to\R$ by
\begin{equation}\label{J}
J(u)= \dfrac{1}{2}\left[\int_0^1 \left(u^2(x) + a(x) u'^2(x)\right)dx + \beta a(1)u^2(1) + a(1) R(u(1)-f(1)) -\int_0^1 (f+g)(x)u(x)dx
\right].
\end{equation}
Then one can check that $J$ is continuously differentiable on $W^1_a(0,1)$ and its differential is given by
\begin{equation}\label{diffJ}
J'(u).\phi=\int_0^1\big(u\phi+au'\phi'\big)dx + \beta a(1) u(1)\phi(1)+ a(1) \rho(u(1)-f(1))\phi(1) - \int_0^1 (f+g) \phi dx \quad \forall \ u, \phi \in
W^1_a(0,1).
\end{equation}
Moreover, since $\rho$ is nondecreasing on $\mathbb{R}$, we deduce that $J$ is a strictly convex function and 
$$
J(u) \geqslant \dfrac{1}{2}|||u|||_{1,a}^2 - ||f+g||_{L^2(0,1)}||u||_{L^2(0,1)} \geqslant |||u|||_{1,a} \left(\dfrac{1}{2}|||u|||_{1,a} - \dfrac{1}{\alpha_a}||f+g||_{L^2(0,1)}\right) \quad \forall \ u \in W^1_a(0,1).
$$
Hence since the norm $|||\cdot|||_{1,a}$ is equivalent to the norm $||\cdot||_{1,a}$ on $W^1_a(0,1)$, $J(u) \longrightarrow +\infty$ as $||u||_{1,a} \longrightarrow +\infty$. Hence $J$ is coercive and strictly convex on $W^1_a(0,1)$ and thus $J$ attains
a minimum at some unique point $u \in   W^1_a(0,1)$, which satisfies the Euler equation
$$
J^{\prime}(u)=0 \,.
$$
Thus $u \in W^1_a(0,1)$ is the unique solution of
\begin{equation}\label{EulerJ}
\int_0^1\big(u\phi+au'\phi'\big)dx + \beta a(1) u(1)\phi(1)+ a(1) \rho(u(1)-f(1))\phi(1) - \int_0^1 (f+g) \phi dx=0 \quad \forall \phi \in
W^1_a(0,1).
\end{equation}
In particular for all $\phi \in \mathcal{C}^{\infty}_c(0,1)$, we have
 $$
 \int_0^1\big(u\phi+au'\phi'\big)dx =\int_0^1(f+g)\phi dx  \qquad \forall \ \phi \in \mathcal{C}^{\infty}_c(0,1) \,.
 $$
 Hence by duality, we have $u-(au')'=f+g$ in the sense of distributions. Thus $u \in  W^2_a(0,1)$ and
 $$
 u-(au')'=f+g \qquad \mbox{a.e in } (0,1)\,.
 $$
 This yields
 $$
 a(1)\phi(1)\left[u'(1) +\beta u(1) + \rho(u(1)-f(1))\right]=0 \qquad \forall \ \phi \in W^1_a(0,1)\,.
 $$
 Since $a(1)>0$ and the function $\phi$ defined by $\phi(x)=x$ for all $x \in (0,1)$ is in $W^1_a(0,1)$ we deduce
 that 
 $$
 u'(1) +\beta u(1) + \rho(u(1)-f(1))=0 \,.
 $$
Setting $v=u-f$, we check that $(u,v) \in D(A_\beta^{nl})$ and solves \eqref{eq:maxNL}.
\end{Proof}
Hence thanks to classical results on nonlinear maximal monotone operators (see e.g. \cite{Brezis, Barbu}), we have
\begin{corollary}
Assume \eqref{eq:hp_a} and that $\rho$ satisfies the above assumptions. Then, for any $U_0=(u_0,u_1)\in D(A_\beta^{nl})$, 
problem~\eqref{DWsNL}-\eqref{b+isNL} has a unique solution $u$ such that
\begin{equation*}
%\label{eq:}
u\in W^{2,\infty}\big([0,\infty);L^2(0,1)\big)\cap W^{1,\infty}\big([0,\infty);W^1_a(0,1)\big)\cap L^{\infty}\big([0,\infty);W^2_a(0,1)\big),
\end{equation*}
Moreover the energy of $u$ defined by \eqref{energy_u}
satisfies the dissipation relation
\begin{equation}\label{dissip_uNL}
\dfrac{d E_u}{dt}(t)=-a(1) u_t(t,1)\rho(u_t(t,1)) \leqslant 0 \qquad \forall \ t \geqslant 0.
\end{equation}
\end{corollary}
\subsection{Nonlinear stability analysis}
We now follow the optimal-weight convexity method introduced in \cite{alaamo2005} and simplified in \cite{alajde2010} (see also~\cite{alaCIME}). For this, we need
to introduce several functions. We first define a function $H: [0,r_0^2] \to [0,\infty)$ by
\begin{equation}\label{H}
H(x) =\sqrt{x} g(\sqrt{x}) \quad x \in [0,r_0^2],
\end{equation}
where $r_0 \leqslant 1$ is assumed to be sufficiently small. We assume that $H$is strictly convex on $[0,r_0^2]$. We extend $H$ to a function $\widehat{H}$ on $[0,\infty)$ by setting $\widehat{H}(x)=+\infty$ when $x \notin [0,r_0^2]$. We then define a function
$L$ on $[0,\infty)$ by
\begin{equation}\label{L}
L(y)=\left\{\begin{array}{ll}
\displaystyle{\frac{\widehat H^\ast(y)}{y}} & \textrm{if }y>0,\\
0 & \textrm{if }y=0,
\end{array}
\right.
\end{equation}
where $\widehat H^{\ast}$ stands for the convex conjugate of $\widehat H$ defined by $\widehat H^\ast (y)=\sup_{x\in \R}\{xy-\widehat H(x)\}$.
One can show that $L$ is a continuous increasing, one-to-one and onto function from $[0,\infty)$ on $[0,r_0^2)$. Moreover $L$ is continuously differentiable on $(0, \infty)$ and 
\begin{equation}\label{gamma}
0 < L(H'(r_0^2))<r_0^2,
\end{equation}
 holds (see \cite{alaamo2005, alajde2010} for more details). Finally we define a function $\Lambda_H$ on $[0,r_0^2]$
\begin{equation}\label{LambdaH}
\Lambda_H(x)=\frac{H(x)}{xH'(x)}.
\end{equation}
Note that $\Lambda_H([0,r_0^2]) \subset [0,1]$ thanks to our convexity assumptions.
\begin{theorem}\label{NLstab}
We assume the above hypotheses on $a$ and on $\rho$, $g$ and $H$, and that $\beta > 0$ is given. Let $(u_0,u_1) \in \mathcal{H}_{\beta}$ be given such that $E_u(0)>0$, and $u$ be the corresponding solution of ~\eqref{DWsNL}-\eqref{b+isNL}. Let $\gamma>
\max( \frac{E_u(0)}{2L(H'(r_0^2))},C_6)$ (where $C_6$ is an explicit constant appearing in~\eqref{stab15NL})
then the energy $E_u$ of $u$ satisfies the following estimate:
\begin{equation}\label{decayenergyE}
\displaystyle{E_u(t)\leqslant 2 \gamma L\Big(\frac{1}{\psi_{  0}^{-1}(\frac{t}{M})}\Big) \ , \quad \forall \ t\geq
\frac{M}{H'(r_0^2)}}\,.
\end{equation}
where
\begin{equation}\label{defpsir}
\psi_0(x)=\frac{1}{H'(r_0^2)}+\int_{1/x}^{H'(r_0^2)}\frac{1}{y^2(1-\Lambda_H((H')^{-1}(\theta)))}\,dy.
\end{equation}
Furthermore, if \, $\limsup_{x \rightarrow 0^+} \Lambda_H(x)<1$, then $E$ satisfies the following simplified decay rate
\begin{equation}\label{decayenergysimpE}
\displaystyle{E_u(t)\leq2 \gamma \Big(H'\Big)^{-1}\Big(\frac{\kappa M}{t}\Big) \ , }
\end{equation}
for $t$ sufficiently large, and where $\kappa >0$ is a constant independent of $E(0)$.
\end{theorem}
\begin{remark}\label{decayexamples}\rm
The above theorem shows that the solutions of the boundary degenerate  nonlinearly damped wave equation above have the same stability properties as the corresponding nondegenerate nonlinearly damped wave equation, that is both have the same decay rates of their energies. In particular, 
\begin{itemize}
\item For the polynomial case for which $g(x)=|x|^{p-1}x$ in a neighborhood of $x=0$ with $p>1$, 
$$
E_u(t) \leqslant C_{E_u(0)} \gamma t^{-\frac{2}{p-1}} \mbox{ for sufficiently large } t.
$$
\item For $g(x)= |x|^{p-1}x \ln^q(\frac{1}{|x|})$ in a neighborhood of $x=0$ with $p>1, q>0$, 
$$
E_u(t) \leqslant C_{E_u(0)} \gamma t^{-\frac{2}{p-1}}(\ln (t))^{-2q/(p-1)} \mbox{ for sufficiently large } t.
$$
\item For $g(x)=sign(x)e^{-1/x^2}$ in a neighborhood of $x=0$, 
$$
E_u(t) \leqslant C_{E_u(0)} \gamma \ln^{-1} (t) \mbox{ for sufficiently large } t.
$$
\item For $g(x)=sign(x) e^{-\ln^p(\frac{1}{|x|})}$ in a neighborhood of $x=0$ with $p>2$
$$
E_u(t) \leqslant C_{E_u(0)} \gamma  e^{-2\left( \ln \left(t\right)\right)^{1/p}}  \mbox{ for sufficiently large } t.
$$
\end{itemize}
Here $\gamma$ is as in Theorem~\ref{NLstab}
 (see e.g. \cite{Komor94} for the linear and polynomial cases and \cite{alajde2010} for the other cases and the references therein, and also \cite{alawanyu} for the last example in the case $p>2$).
\end{remark}
%%%%%%%%%%%%%%%%%%%%
%%%%%%%%%%
 \begin{Proof}
 Thanks to the density of $D(A_\beta^{nl})$ in $\mathcal{H}_{\beta}$, and since $A_\beta^{nl}$ is a maximal dissipative operator, it is sufficient to consider smooth initial data $(u_0,u_1)$. Hence, let $U_0=(u_0,u_1)\in D(A_\beta^{nl})$ be given, and $u$ be the corresponding solution of problem~\eqref{DWsNL}-\eqref{b+isNL}. Let $\gamma>
 \frac{E_u(0)}{2L(H'(r_0^2))}$ which will be precise later on in the proof and define the optimal-weight function as
 \begin{equation}\label{optw}
 w(s)= L^{-1}\left(\frac{E_u(s)}{2\gamma}\right) \quad \forall \ s \geqslant 0.
 \end{equation}
 We multiply \eqref{DWsNL} by $w(E_u(t))xu_x$ and integrate the resulting equation over $(S,T) \times (0,1)$. After suitable integrations by parts as in the previous section, this gives  for all $0 \leqslant S \leqslant T$
 \begin{multline*}
 \int_S^Tw(E_u(t))\int_0^1 \Big(-x \big(\dfrac{u_t^2}{2}\big)_x + a(x) u_x^2 + xa(x) \big(\dfrac{u_x^2}{2}\big)_x\Big)dxdt + \Big[w(E_u(t))\int_0^1xu_xu_t dx
 \Big]_S^T -\int_S^T w(E_u(t))\big[xa u_x^2\big]_0^1 dt -\\
\int_S^Tw'(E_u(t))E_u'(t) \int_0^1 xu_xu_t dx dt =0.
 \end{multline*}
 We integrate by parts twice again. This gives, together with the trace results at $x=0$ as in the previous section
 \begin{multline}\label{stab1NL}
 \int_S^Tw(E_u(t))\int_0^1 \Big(\dfrac{u_t^2}{2} + (a-x a') \dfrac{u_x^2}{2}\Big)dxdt + \Big[w(E_u(t))\int_0^1xu_xu_t dx
 \Big]_S^T -\dfrac{1}{2} \int_S^T w(E_u(t))\Big(a(1) u_x^2(t,1) + u_t^2(t,1)\Big)dt -\\
 \int_S^Tw'(E_u(t))E_u'(t) \int_0^1 xu_xu_t dx dt=0 \quad \forall \ 0 \leqslant S \leqslant T.
 \end{multline}
 We multiply \eqref{DWsNL} by $w(E_u(t))u$ and integrate the resulting equation over $(S,T) \times (0,1)$. This gives after a suitable integration by parts and thanks to our trace results.
 \begin{multline}\label{stab2NL}
 \int_S^Tw(E_u(t))\int_0^1 \Big(-u_t^2 + au_x^2\Big) dxdt + \Big[w(E_u(t))\int_0^1u_t udx
 \Big]_S^T - \int_S^T w(E_u(t))a(1)u_x(t,1)u(t,1)dt-\\
\int_S^Tw'(E_u(t))E_u'(t) \int_0^1uu_t dx dt  =0 \quad \forall \ 0 \leqslant S \leqslant T.
 \end{multline}
 We now combine \eqref{stab1NL} multiplied by 2 with \eqref{stab2NL} multiplied $\dfrac{\mu_a}{2}$. This gives
 \begin{multline}\label{stab3NL}
 \int_S^Tw(E_u(t))\int_0^1 \Big[(2-\mu_a)\dfrac{u_t^2}{2} + \big[2(a-x a')+ a\mu_a\big] \dfrac{u_x^2}{2}\Big]dxdt + \dfrac{2-\mu_a}{2}\beta a(1)\int_S^T w(E_u(t))u^2(t,1)dt=\\
 -2 \Big[w(E_u(t))\int_0^1xu_xu_t dx \Big]_S^T - \dfrac{\mu_a}{2}\Big[w(E_u(t))\int_0^1u_tu dx \Big]_S^T+\\
 \int_S^T w'(E_u(t))E_u'(t) \int_0^1\left(2xu_x + \dfrac{\mu_a}{2} u\right)u_t dxdt +
 \int_S^Tw(E_u(t)) \tilde{h}(t)dt \quad \forall \ 0 \leqslant S \leqslant T,
 \end{multline}
 where the function $\tilde{h}$ is given by
 \begin{equation}\label{defhNL}
 \tilde{h}(t)=u_t^2(t,1) + a(1) \rho(u_t(t,1))^2 + a(1) \beta (1+\beta-\mu_a)u^2(t,1) + \big(2\beta-\dfrac{\mu_a}{2}\big)a(1) \rho(u_t(t,1))u(t,1) \quad t \in (S,T).
 \end{equation}
 By definition of $\mu_a$, we have
 $$
 (2-\mu_a)a \leqslant 2(a-x a')+ a\mu_a.
 $$
 This, together with \eqref{stab3NL}, gives
 \begin{multline}\label{stab4NL}
(2-\mu_a) \int_S^TE_u(t)dt \leqslant 
 - \Big[w(E_u(t))\int_0^12xu_xu_t  + \dfrac{\mu_a}{2} u_tu dx \Big]_S^T + \int_S^T w'(E_u(t))E_u'(t) \int_0^1\left(2xu_x + \dfrac{\mu_a}{2} u\right)u_t dxdt +\\
 \int_S^T \tilde{h}(t)dt \quad \forall \ 0 \leqslant S \leqslant T.
 \end{multline}
On the other hand, we have
\begin{equation}\label{stab5NL}
\tilde{h}(t) \leqslant \eta_3 u_t^2(t,1) + \eta_4\rho(u_t(t,1))^2 + \eta_5 a(1)u^2(t,1) \quad \forall \  t \in (S,T),
\end{equation}
where $\eta_i$ for $i=3,4,5$ are positive constants which do not depend on the weight function $w$ nor on $E(t)$.
The two first terms in \eqref{stab4NL} are estimates as in the linear stabilization case. This, together with the properties that $w$ is nondecreasing whereas $E$ is non increasing yield
\begin{multline}\label{stab6NL}
(2-\mu_a) \int_S^TE_u(t)dt \leqslant K_aw(E_u(S))E_u(S)+ \eta_3\int_S^T w(E_u(t))u_t^2(t,1)dt + \eta_4 \int_S^Tw(E_u(t))\rho(u_t(t,1))^2dt +\\
\eta_5 \int_S^Tw(E_u(t)) a(1) u^2(t,1)dt. 
\end{multline}
where $K_a$ is a positive constant which do not depend on the weight function $w$ nor on $E(t)$.
We now estimate the last term of this inequality as in the linear stabilization case, using once again in addition our optimal weight function. Set $\lambda=u(t,1)$ and denote by $z$ the solution of the
degenerate elliptic problem \eqref{eq: var0s}. We multiply \eqref{DWsNL} by $w(E_u(t))z$ and integrate the resulting equation over $(S,T) \times (0,1)$. This gives after suitable integrations by parts.
\begin{multline}\label{stab8NL}
\int_S^Ta(1)w(E_u(t))u^2(t,1)dt=\int_S^Tw(E_u(t))\int_0^1 u_tz_t dxdt+\int_S^Tw'(E_u(t))E_u'(t)\int_0^1 u_tzdxdt -\\
a(1)\int_S^T w(E_u(t))\rho(u_t(t,1))z(t,1)dt-\Big[w(E_u(t))\int_0^1u_t zdx\Big]_S^T.
\end{multline}
We now estimate the terms of the right hand side in this inequality, as follows.
Using \eqref{stab9}-\eqref{stab11} in \eqref{stab8NL}, we obtain for all $\delta>0$
\begin{multline*}
\int_S^Ta(1)w(E_u(t))u^2(t,1)dt \leqslant \delta \int_S^T w(E_u(t))E_u(t)dt +C_1w(E_u(S))E_u(S) +  \\
C_2 (1+ \dfrac{1}{\delta})\int_S^T w(E_u(t))\left(\rho(u_t(t,1))^2+ u_t^2(t,1)\right)dt,
\end{multline*}
where $C_1, C_2$ are positive constants which do not depend on the weight function $w$ nor on $E(t)$.
Choosing $\delta =\dfrac{2-\mu_a}{2\eta_5}$ in the above inequality and combining the resulting
inequality in \eqref{stab6NL} yield
\begin{equation}\label{stab12NL}
\int_S^Tw(E_u(t))E_u(t)dt \leqslant  C_3 w(E_u(S))E_u(S) + C_4 \int_S^Tw(E_u(t))\left(\rho^2(u_t(t,1))+ u_t^2(t,1)\right)dt,
\end{equation}
 where $C_3, C_4$ are positive constants which do not depend on the weight function $w$ nor on $E(t)$.
 It remains to estimate the last term on the right hand side of the above inequality. We further proceed as in \cite{alaamo2005, alajde2010}. That is 
 we fix $t \geqslant 0$. Assume first that $|u_t(t,1)| \leqslant \varepsilon_0$ where $\varepsilon_0=\min(1,g(r_0))$.
 Hence, thanks to our assumption on $\rho$, we have 
 $$
 \left|\dfrac {|\rho(u_t(t,1)|}{c_2}\right|^2 \leqslant \left|g^{-1}(u_t(t,1))\right|^2 \leqslant |g^{-1}(\varepsilon_0)|^2 \leqslant r_0^2.
 $$
 On the other hand, we have
 $$
 H\left(\dfrac{|\rho(u_t(t,1)|^2}{c_2^2}\right)=\dfrac{|\rho(u_t(t,1)|}{c_2}g\left(\dfrac{|\rho(u_t(t,1)|}{c_2}\right) \leqslant \dfrac{1}{c_2}u_t(t,1)\rho(u_t(t,1)).
 $$
 Hence, since $H$ is nondecreasing, we have whenever $t$ is such that $|u_t(t,1)| \leqslant \varepsilon_0$
 \begin{equation}\label{stab13NL}
w(E_u(t)) |\rho(u_t(t,1)|^2 \leqslant c_2^2w(E_u(t))H^{-1}\left(\dfrac{1}{c_2}u_t(t,1)\rho(u_t(t,1)) \right) \leqslant
c_2^2 \widehat{H}^{\ast}(w(E_u(t))) + c_2 u_t(t,1)\rho(u_t(t,1)).
\end{equation}
 We now assume that $t$ is such that $|u_t(t,1)| \geqslant \varepsilon_0$, then up to a change in the constants $c_1$ and $c_2$ in \eqref{hyprho}, we can assume
 $$
 |\rho(u_t(t,1))|\leqslant c_2 |u_t(t,1)|,
 $$
so that
$$
\int_{ t \in [S,T], |u_t(t,1)| \geqslant \varepsilon_0} w(E_u(t))|\rho(u_t(t,1))|^2 \leqslant \dfrac{c_2}{a(1)}   w(E_u(S)) E_u(S). 
$$ 
Combining this last estimate together with \eqref{stab13NL}, we obtain
\begin{equation}\label{stab14NL}
\int_S^Tw(E_u(t)) |\rho(u_t(t,1))|^2dt \leqslant c_2^2 \int_S^T \widehat{H}^{\ast}(w(E_u(t))) + \dfrac{c_2}{a(1)}E_u(S)\left(1+w(E_u(S))\right).
\end{equation}
We similarly estimate the term $\int_S^Tw(E_u(t))u_t^2(t,1)dt$ proceeding as in \cite{alaamo2005, alajde2010}. That is, we fix $t \geqslant 0$. We consider first the case for which $|u_t(t,1)| \leqslant \varepsilon_1$ where $\varepsilon_{1}=\min \{r_{0},g(r_{1})\}$ where $r_{1}$ is defined by
$$
r_{1}^2=H^{-1}\left(\frac{c_1}{c_2}H(r_{0}^2)\right).
$$
Thanks to our assumptions on $\rho$, we have
$$
 H\left(|u_t(t,1)|^2\right) \leqslant \dfrac{1}{c_1} u_t(t,1)\rho(u_t(t,1).
$$
Hence, we have
\begin{equation}\label{stab15NL}
w(E_u(t)) |u_t(t,1)|^2 \leqslant w(E_u(t))H^{-1}\left(\dfrac{1}{c_1}u_t(t,1)\rho(u_t(t,1)) \right) \leqslant
\widehat{H}^{\ast}(w(E_u(t))) + \dfrac{1}{c_1} u_t(t,1)\rho(u_t(t,1)).
\end{equation}
Assume now that $t$ is such that $|u_t(t,1)| \geqslant \varepsilon_1$, then up to a change in the constants $c_1$ and $c_2$ in \eqref{hyprho}, we can assume
 $$
 |\rho(u_t(t,1))|\geqslant c_1 |u_t(t,1)|,
 $$
so that
$$
\int_{ t \in [S,T], |u_t(t,1)| \geqslant \varepsilon_1} w(E_u(t))|u_t(t,1)|^2 \leqslant \dfrac{1}{c_1a(1)}   w(E_u(S)) E_u(S). 
$$ 
Combining this last estimate together with \eqref{stab13NL}, we obtain
\begin{equation}\label{stab16NL}
\int_S^Tw(E_u(t)) |u_t(t,1)|^2dt \leqslant  \int_S^T \widehat{H}^{\ast}(w(E_u(t)))dt + \dfrac{1}{c_1a(1)}E_u(S)\left(1+w(E_u(S))\right).
\end{equation}
On the other hand, we recall that $\gamma$ satisfies
\eqref{gamma}, thus we have
$$
w(E_u(S)) \leqslant L^{-1}\left(\dfrac{E_u(0)}{2\gamma}\right) <H'(r_0^2) \quad \forall \ S \geqslant 0.
$$
Inserting the estimates \eqref{stab14NL} and \eqref{stab16NL} in \eqref{stab12NL}, and using the above estimate, we obtain
\begin{equation}\label{stab17NL}
\int_S^Tw(E_u(t))E_u(t)dt \leqslant  C_5 E_u(S) + C_6\int_S^T \widehat{H}^{\ast}(w(E_u(t)))dt,
\end{equation}
where $C_5, C_6$ are positive constants which do not depend on the weight function $w$ nor on $E(t)$. Thanks to our choice of weight function w
$$
L(w(E_u(t)))=\dfrac{E_u(t)}{2\gamma} \quad \forall  \ t \geqslant 0,
$$
so that we have
$$
\int_S^Tw(E_u(t))E_u(t)dt \leqslant  C_5 E_u(S) + \dfrac{C_6}{2\gamma}\int_S^T w(E_u(t))E_u(t)dt,
$$
Choosing $\gamma \geqslant C_6$ in addition to \eqref{gamma}, we obtain that
\begin{equation}\label{stab18NL}
\int_S^Tw(E_u(t))E_u(t)dt \leqslant  M E_u(S) \quad \forall \ 0 \leqslant S \leqslant T.,
\end{equation}
where $M=2C_5$.
 Then proof can be completed applying the following result (see \cite[Theorem 2.3] {alajde2010}). \end{Proof}
 \begin{theorem}\label{thmintc} 
 %\cite[Theorem 2.3]{alajde2010}
Let $H$ be a strictly convex function on $[0,r_0^2]$ such that $H(0)=H'(0)=0$ and define $L$ and $\Lambda_H$ as above.
Let $E$ be a given nonincreasing, absolutely continuous
 function from $[0,+\infty)$ on $[0,+\infty)$ with $E(0)>0$ satisfying the
 following weighted nonlinear inequality
\begin{equation}\label{ineqint1}
\int_{S}^T L^{-1}(\frac{E(t)}{2\gamma}) E(t)  \, dt \leq M E(S) \,, \quad \forall \, 0\leq S \leq T.
\end{equation}
where $M>0$ and where $\gamma> \dfrac{E(0)}{2L(H'(r_0^2))}$. 
Then $E$ satisfies the following estimate:
\begin{equation}\label{decayenergy}
\displaystyle{E(t)\leqslant 2 \gamma L\Big(\frac{1}{\psi_{  0}^{-1}(\frac{t}{M})}\Big) \ , \quad \forall \ t\geq
\frac{M}{H'(r_0^2)}}\,.
\end{equation}
where $\psi_0$ is defined in \eqref{defpsir}.
Furthermore, if \, $\limsup_{x \rightarrow 0^+} \Lambda_H(x)<1$, then $E$ satisfies the following simplified decay rate
\begin{equation}\label{decayenergysimp}
\displaystyle{E(t)\leq2 \gamma \Big(H'\Big)^{-1}\Big(\frac{\kappa M}{t}\Big) \ , }
\end{equation}
for $t$ sufficiently large, and where $\kappa >0$ is a constant independent of $E(0)$.
\end{theorem}

\begin{remark} \label{re:Lyapunov}\rm
It should be noted that one can also reformulate, with no mathematical originality and no gain with respect to applications and research, all our results on the nonlinear stabilization of degenerate equations of this section by means of a "Lyapunov" presentation. In this case,
it is sufficient to track all the steps of our proof, remove all the integrations with respect to time (from $S$ to $T$) and  multiply afterwards the resulting inequality by a weight function, which can be a weaker (and less good) weight function than in the original method introduced for the first time in \cite{alaamo2005} (see also \cite{alajde2010}). This weaker weight function can easily be deduced by dropping in the original
computations of \cite{alaamo2005}, the negative part in the convex conjugate of the strictly convex function $H^{\ast}$ defined in \eqref{H}. Namely,
this consists in replacing $H^{\ast}(y)=y(H^{\prime})^{-1}(y) -
H((H^{\prime})^{-1}(y))$ for $ y \in [0, c]$ (for a suitable $c>0$) in the original paper by the function $H_2(y)=y (H^{\prime})^{-1}(y)$ for $ y \in [0, c]$.
 The results would also be weaker and destroy some nice and further properties proved later on in \cite{alajde2010} which lead to simplified and {\em optimal} energy decay rates.
\end{remark}

%%%%%%%%%%%%%
%%%%%%%%%%%%%%%%%%%%

\bibliographystyle{plain}
\bibliography{biblio_DW}

\end{document}